\documentclass[preprint,1p,11pt]{IR-Template/ISAS_IR}
\usepackage{amssymb, natbib}
\usepackage{amsmath, amsthm}

\usepackage{graphicx}
\usepackage{listings}
\usepackage{color}
\usepackage[TABBOTCAP]{subfigure}

\usepackage[nesting]{hyperref}

\newcommand{\BlackBoxes}{\global\overfullrule5pt}

\BlackBoxes

\newcommand{\R}{\mathbb{R}}
\newcommand{\N}{\mathbb{N}}

\newcommand{\PP}{\mathbb{P}}

\newcommand{\Pop}{\operatorname{\mathbb{\PP}}}

\newcommand{\FF}{\mathcal{F}}

\newtheorem{theorem}{Theorem}

\theoremstyle{definition}

\numberwithin{equation}{section} \numberwithin{theorem}{section}

\def\0{\kern0pt\-\nobreak\hskip0pt\relax}

\makeatletter
\AtBeginDocument{%
 \def\@serieslogo{%
 \vbox to\headheight{%
 \parindent\z@ \fontsize{6}{7\p@}\selectfont
\today\; Final 
 \vss}}}

\def\makeoverbar#1#2#3#4#5#6#7{%
 \setbox0=\hbox{$\m@th#2\mkern#5mu{{}#3{}}\mkern#6mu$}%
 \setbox1=\null \dimen@=#4\fontdimen8#13 \dimen@=3.5\dimen@
 \advance\dimen@ by \ht0 \dimen@=-#7\dimen@ \advance\dimen@ by \wd0
 \ht1=\ht0 \dp1=\dp0 \wd1=\dimen@
 \dimen@=\fontdimen8#13 \fontdimen8#13=#4\fontdimen8#13
 \rlap{\hbox to \wd0{$\m@th\hss#2{\overline{\box1}}\mkern#5mu$}}
 \fontdimen8#13=\dimen@}

\def\mylabel#1#2{{\def\@currentlabel{#2}\label{#1}}}

\makeatother

\begin{document}


\makeatletter \providecommand\@dotsep{5} \makeatother

\title{Exact and Approximate Hidden Markov Chain Filters Based on Discrete Observations}

\author[stoch]{Nicole B\"auerle}
\ead{nicole.baeuerle@kit.edu}
\address[stoch]{Institute of Stochastics,
Karlsruhe Institute of Technology, D-76131 Karlsruhe, Germany}

\author[isas]{Igor Gilitschenski}
\ead{gilitschenski@kit.edu}

\author[isas]{Uwe D. Hanebeck}
\ead{uwe.hanebeck@ieee.org}
\address[isas]{Institute for Anthropomatics and Robotics,
Karlsruhe Institute of Technology, D-76131 Karlsruhe, Germany}

\begin{abstract}
We consider a {\em Hidden Markov Model (HMM)} where the integrated continuous-time Markov chain can be observed at discrete time points perturbed by a Brownian motion. The aim is to derive a filter for the underlying continuous-time Markov chain. The recursion formula for the discrete-time filter is easy to derive, however involves densities which are very hard to obtain. In this paper we derive exact formulas for the necessary densities in the case the state space of the HMM consists of two elements only. This is done by relating the underlying integrated continuous-time Markov chain to the so-called asymmetric telegraph process and by using recent results on this process. In case the state space consists of more than two elements we present three different ways to approximate the densities for the filter. The first approach is based on the continuous filter problem. The second approach is to derive a PDE for the densities and solve it numerically and the third approach is a crude discrete time approximation of the Markov chain. All three approaches are compared in a numerical study.
\end{abstract}


\maketitle

\section{Introduction}
We consider a {\em Hidden Markov Model (HMM)} where the integrated continuous-time Markov chain can be observed at discrete time points perturbed by a Brownian motion. The aim is to derive a filter for the underlying continuous-time Markov chain. Thus, we have a continuous-time Hidden Markov Model with discrete observations.  Models of this type are widely used in finance, in telecommunication or in biology. For references, we refer to the examples in Section \ref{sec:filter1}. The knowledge of the filter is for example necessary in case one has to control the system in a dynamic way. Usually the principle of estimation and control is valid and implies that the filter has to be computed first and then the control is applied, see e.g. \cite{BL09,br,SH04} for different kinds of applications.
We assume that all parameters of the model are known or have been estimated before. In particular we assume that the number and values of the states of the hidden Markov chain are known. As far as model selection is concerned we refer e.g.\ to \cite{FS01,OG02,DH08}. Papers dealing with the estimation of state values and transition intensities are e.g. \cite{ryden96,EKS08,HS09,HFSS10}.  Since we assume a finite number of states for the Markov chain, the filter has a finite dimension and it is in principal clear form the existing theory how to construct it. However, the filter recursion involves densities which are difficult to compute explicitly.

We rely on results from \cite{LR14} and \cite{MKR94} on the asymmetric telegraph process to obtain the exact filter in closed form in case the hidden Markov chain has only two states. We distinguish between the symmetric case, where the switching intensities are the same for both states and the asymmetric case where intensities are different. The first case is easier and related to the well-known {\em telegraph process}. The second problem is related to the asymmetric telegraph process. Estimates for the switching intensity of a symmetric telegraph process have e.g. been derived in \cite{Yao85}, \cite{IY08} and for an inhomogeneous telegraph process in \cite{I01}.

In the general situation with more than two states, the required densities become rather complicated. So a procedure to compute these densities approximately is called for. We present three different ways to derive approximations. All three approaches are general and need no information about the size of the intensities, i.e. how fast the hidden Markov chain switches. Of course the faster the switching of the  Markov chain, the harder it is to obtain a reasonable filter, given a fixed time lag between observations. If the switching intensity is too large compared to the observation time lag, the stationary distribution is a reasonable approximation. The first approach relies on results in \cite{PR10,KPS93}  where the authors derive an approximate continuous-time filter. When this continuous-time filter is discretized, it can also be used by someone who is only able to observe the data in discrete time. However, this approach is based on an assumption which is not true in general. The second approach is to derive a PDE for the densities which are involved in the filter process and to solve this PDE numerically. This approach is exact up to the error of the numerical PDE solution. The third attempt relies on a naive approximation of the continuous-time Markov chain by a discrete-time Markov chain. The finer the approximation, i.e. the smaller the time step of the discrete-time Markov chain, the better the approximate filter but the higher the computational effort.

Our paper is organized as follows: The next section contains a precise mathematical formulation of the problem and derives the basic filter recursion. Section
\ref{sec:2states} is then devoted to the problem with two states. Results about the telegraph process are used to derive the filter explicitly in both the symmetric and the asymmetric case. The following section considers the general problem and three approximate filters are derived. The last section is devoted to numerical results. In particular we compare our approximate filters in a setting with five states in situations with frequent or rare observations and with small or large variability of the error terms.

\section{Hidden Markov Chain Filters}\label{sec:filter1}
Suppose we have a filtered probability space $(\Omega, \FF_T, (\FF_t),\Pop)$ with $T>0$ and a hidden stationary continuous-time Markov chain $\varepsilon=(\varepsilon_t)_{t\in[0,T]}$ on it with finite state space $\mathcal{E}=\{1,\ldots ,d\}$ and $d\in\N$. The intensity matrix of this process is given by $Q$ and the initial distribution by $p_0$. Let us denote by $p(t)=(p_1(t),\ldots,p_d(t))^\top$ the distribution of $\varepsilon_t$ for $t\in[0,T]$, then $p(t)$ satisfies the ODE \begin{equation}\label{eq:Q}\frac{ d p(t)}{dt } = Q p(t)\end{equation} with $p(0)=p_0$. The formal solution of \eqref{eq:Q} is given by  \begin{equation*}p(t) = e^{Qt} p_0 \end{equation*} where  \begin{equation*} e^{Qt} = \sum_{k=0}^\infty \frac{(Qt)^k }{k!} \end{equation*} is the matrix exponential. Moreover, the $h$-step transition probabilities are given by  \begin{equation}\label{eq:transprob} p_{ij}(h) = \Pop(\varepsilon(t+h) = j | \varepsilon(t)=i) = (e^{Qh})_{i,j}. \end{equation}
Now assume further that $W=(W_t)_{t\in[0,T]}$ is a Brownian motion on our filtered probability space which is independent of the hidden continuous-time Markov chain $\varepsilon$ and plays the role of a noise. Moreover, $\alpha_1,\ldots ,\alpha_d$ are different real numbers and $\sigma>0$. We assume that we can observe the process
\begin{equation}\label{eq:Z} Z_t := \int_0^t \alpha_{\varepsilon_s}ds + \sigma W_t \end{equation}
at discrete time points $t_k = kh$ for $k=0,1,\ldots ,N$, $h>0$ and $Nh = T$, i.e. the observation $\sigma$-algebra at time $t_k$ is given by  \begin{equation*}\mathcal{Z}_k := \sigma(\{ Z_{t_1}-Z_{t_0},\ldots ,Z_{t_k}-Z_{t_{k-1}}\}),\quad k=1,\ldots,N \end{equation*} as opposed to the filtration $(\FF_t)$ which contains both the information about $\varepsilon$ and $W$.
A simple interpretation of this model is as follows: Suppose a particle is moving on the real line, starting in zero. The velocity of the particle at time $t$ is determined by $\varepsilon_t$. Thus, $\int_0^t \alpha_{\varepsilon(s)}ds $ gives the exact position of the particle at time $t$. However, we are only able to observe the position at discrete time points $t_k$ together with a noise $\sigma W_{t_k}$. The aim is now to filter from this observation the current velocity of the particle, i.e. to determine
 \begin{equation*} \Pop(\varepsilon_{t_n} = j | \mathcal{Z}_n) =: \mu_n(j) \end{equation*}
where $\mu_0=p_0$ is given. It is well-known that we can derive a recursive filter formula with the help of Bayesian analysis which is given by
\begin{equation}\label{filter} \mu_n(j) = \frac{\sum_{i=1}^d p_{ij}(h) g_{ij}(z,h) \mu_{n-1}(i)}{\sum_{i=1}^d g_{i}(z,h) \mu_{n-1}(i)}\end{equation}
where $g_i(z,h)$ is the density of the conditional distribution of $Z_{h}$ given $\varepsilon_0=i$ and
$g_{ij}(z,h)$ is the density of the conditional distribution of $Z_h$ given $\varepsilon_0=i$ and $\varepsilon_h=j$, see e.g. \cite{EAM95}, chapter 2, \cite{frikr07} chapter 3 and \cite{br11} chapter 5, i.e.
\begin{eqnarray*}
  g_i(z,h) &=& \Pop_i(Z_h \in dz)/dz, \quad i\in \mathcal{E} \\
  g_{ij}(z,h) &=& \Pop(Z_h \in dz | \varepsilon_0=i, \varepsilon_h=j)/dz,\quad i,j\in \mathcal{E}.
\end{eqnarray*}
Obviously by conditioning we have the relation
\begin{equation}\label{eq:gifromgij}
g_i(z,h) = \sum_{j\in \mathcal{E}} g_{ij}(z,h) p_{ij}(h).
\end{equation}
Thus, it is in principle enough to derive the densities $g_{ij}$  from which we obtain $g_i$.
 In case $\mathcal{E}=\{1,2\}$ consists of two elements only, we will derive explicit expressions for $g_i$ and $g_{ij}$ in the next section. For general $\mathcal{E}$ we derive approximate expressions in section \ref{sec:approx}.\\

Let us also shortly mention the situation when we have continuous observations. In this case the observation filtration is given by
\begin{equation*}\mathcal{Z}_t := \sigma(\{ Z_{s}, s\in[0,t]\}),\quad t\in[0,T] \end{equation*}
and we have to determine
\begin{equation*}\Pop(\varepsilon_{t} = j | \mathcal{Z}_t),\quad t\in[0,T]. \end{equation*}
This is the well-known {\em Wonham filter problem}. From the {\em Kallianpur-Striebel formula} we have the following representation
\begin{equation}\label{eq:KS} \Pop(\varepsilon_{t} = j | \mathcal{Z}_t) = \frac{\xi_j(t)}{\sum_{i=1}^d \xi_i(t)}
\end{equation}
where $(\xi_t)_{t\in[0,T]} := (\xi_1(t),\ldots ,\xi_d(t))_{t\in[0,T]}$ satisfies the {\em Zakai equations}
\begin{equation}\label{eq:zakai}
\xi_t = p_0 + \int_0^t Q\xi_s ds + \int_0^t D \xi_s dZ_s
\end{equation}
for $t\in[0,T]$ where $D$ is a diagonal matrix with elements ${\alpha_1}/{\sigma^2},\ldots {\alpha_d}/{\sigma^2}$ on the diagonal (see e.g. \cite{PR10}). Obviously, this filter makes use of a continuous observation of $(Z_t)$. However, we will later see that approximations of this filter at time t may only depend on $Z_t$ which allows us to use this filter in case of discrete observations, too.

Applications of this model are given next.\\[0.2cm]

{\bf Example 1:} Financial data suggest that parameter of asset prices depend on external macroeconomic factors which may be described by a continuous-time Markov chain (cp. \cite{RTA}). A popular model for example is to take the classical Black Scholes model with parameters driven by a factor processes which is represented by a continuous-time Markov chain. Thus, the log-return of an asset price $S_t$ at time $t$ would be
\begin{equation*}
\log S_t = \int_0^t \Big(\alpha_{\varepsilon_s} -\frac12 \sigma^2\Big) ds + \sigma W_t
\end{equation*}
where $(W_t)$ is a Brownian motion independent of $(\varepsilon_t)$, $\alpha$ is the drift or appreciation rate and $\sigma$ is a fixed volatility. The parameter $\sigma$ may also depend on $(\varepsilon_t)$. Models where the factor process is assumed to be known are among others treated in \cite{MKR94,br04,YZ04,Costa08,ZSM10}. If $(\varepsilon_t)$ is not observable and $\sigma$ constant we have exactly the generic situation which we consider in this paper. The model is then reasonably described by an HMM or Markov switching models for price parameters. The aim is to filter the underlying economic factor. For applications of this model see \cite{SH04,br,HS09,sw10,buv12,fgw12}. In this case it is often only possible to observe the log-return, i.e. the integrated HMM  at discrete time points. For example some asset prices are only quoted on a daily basis.\\[0.4cm]

{\bf Example 2:} In communication networks, information arrives to a multiplexer, switch or information processor at a rate which changes randomly and often shows a high degree of correlation in time. For example \cite{hl86} modeled the input stream of a statistical multiplexer consisting of a mix of data and packetized voice sources as a Markov-modulated Poisson process. In \cite{se91} the authors used a Markov-modulated continuous flow model for the information stream at a multiplexer and analyzed the system performance. These and other findings fostered the investigation of so-called fluid queues, see e.g. \cite{k97}. A more recent reference is \cite{rs04} where the authors are interested in the stationary queueing level of a fluid queue whose buffer content evolution process $(Q_t)$  is given by the SDE
\begin{equation*}
    d Q_t = \alpha_{\varepsilon_t} dt - \mu_{\varepsilon_t} Q_t dt + \sigma_{\varepsilon_t} Q_t d W_t.
\end{equation*}
Here again, if one is able to observe, perturbed by an error, the amount of fluid which came in, then the filter computes the distribution of the current inflow rate.
Such fluid queues can also be used to model production processes.

\section{HMM problem with two states}\label{sec:2states}
\subsection{Symmetric $Q$-matrix}
In case the state space consists of two elements only, i.e.  $\mathcal{E}=\{1,2\}$ and the intensities are symmetric, i.e.
\begin{equation}
Q = \left ( \begin{array}{cc} -\lambda & \lambda\\ \lambda & -\lambda\end{array}\right)
\end{equation}
with $\lambda >0$, the model can be reduced to the so-called {\em telegraph process}. The telegraph process $(I_t)_{t\in[0,T]}$ is given by
\begin{equation}\label{telegraph}
I_t := \int_0^t (-1)^{N_s} ds
\end{equation}
where $(N_t)_{t\in[0,T]}$ is a Poisson process with intensity $\lambda$. The telegraph process has been introduced by Kac in lecture notes from 1956, see also \cite{Kac}. It has been generalized in different ways and studied in various applications. For a recent paper with many references see \cite{LR14}. One usually distinguishes the {\em symmetric} case where the particle moves with symmetric velocities $+1$ and $-1$ and switches the direction with same intensity $\lambda>0$ and the {\em asymmetric} case where we have different arbitrary velocities and different switching intensities.

Now let us denote
\begin{equation}\label{eq:Jt} J_t := \int_0^t \alpha_{\varepsilon_s}ds \end{equation}
and define $a := \frac12 (\alpha_1+\alpha_2), b_1 = \frac12 (\alpha_1-\alpha_2)$ and $b_2 = \frac12 (\alpha_2-\alpha_1)$. Then, given $\varepsilon_0= i$ for $ i=1,2$ we have
\begin{equation}
J_t = a t + b_i I_t,\quad t\in[0,T],
\end{equation}
i.e. $J_t$ is linear transformation of $I_t$. The transition probabilities in \eqref{eq:transprob} in this case simplify to
\begin{equation}
  p_{12}(t) = p_{21}(t) = \frac12 - \frac12 e^{-2\lambda t},\quad t\in[0,T].
\end{equation}
It is now possible to derive explicit formulas for the densities  $g_i, g_{ij}, i,j\in\{1,2\}$. which appear in the filter formula \eqref{filter}.
The density for the telegraph process at time $h$, $I_h$ given the initial state can already be found in the appendix of \cite{MKR94}. In order to derive the densities for our application, let us introduce the following notations. The modified Bessel functions are given by
\begin{equation}\label{eq:Bessel} B_0(z)= \sum_{k=0}^\infty \frac{(z^2/4)^k}{k!^2},\quad  B_1(z)= \frac z 2 \sum_{k=0}^\infty \frac{(z^2/4)^k}{k! (k+1)!},\quad z\ge 0
\end{equation}
and we denote the density of the normal distribution with zero expectation and standard deviation $\sigma$ by
\begin{equation*}
\phi_\sigma(x) = \frac{1}{\sqrt{2\pi}\sigma} e^{-\frac12 \frac{x^2}{\sigma^2}}, \quad x\in\R.
\end{equation*}
Then we obtain:

\begin{theorem}\label{theo:Qsym}
The density of $Z_h$ given $\varepsilon_0=i$ for $i=1,2$ is
\begin{equation*}\label{eq:densityg_isym}
  g_i(z,h) =  e^{-\lambda h} \Big[\phi_{\sigma\sqrt{h}}(z-h(b_i+a)) + \frac{1}{|b_i|} \int_{h (a- |b_i|)}^{h (a+|b_i|)}\phi_{\sigma\sqrt{h}}(z-x) \tilde{q}_{h,\lambda}\Big(\frac{x-ah}{b_i}\Big) dx\Big], \quad z\in\R,
\end{equation*}
with
\begin{equation*}\label{eq:htl} \tilde{q}_{h,\lambda}(x) = \frac\lambda 2  \Big[B_0\Big( \lambda \sqrt{(h^2-x^2)}\Big) + \Big(\frac{h+x}{h-x} \Big)^{\frac 12} B_1\Big( \lambda \sqrt{(h^2-x^2)}\Big)  \Big], \quad x\in [-h,h].\end{equation*}
The density of $Z_h$ given $\varepsilon_0=\varepsilon_h=i$ is
\begin{equation*}\label{eq:densityg_iisym}
  g_{ii}(z,h) = \frac{1}{\cosh(\lambda h)} \Big[\phi_{\sigma\sqrt{h}}(z-h(b_i+a)) + \frac{1}{|b_i|} \int_{h (a-|b_i|)}^{h (a+|b_i|)}\phi_{\sigma\sqrt{h}}(z-x) \tilde{q}_{h,+}\Big(\frac{x-ah}{b_i}\Big) dx\Big], \quad z\in\R,
\end{equation*}
with
\begin{equation*}
  \tilde{q}_{h,+}(x) =  \frac{\lambda}2 \Big(\frac{h+x}{h-x} \Big)^{\frac 12} B_1\Big(\sqrt{(h^2-x^2)} \lambda \Big), \quad x\in [-h,h].
\end{equation*}
The density of $Z_h$ given $\varepsilon_0=i, \varepsilon_h=j$ with $i,j\in\{1,2\}$ and $i\neq j$ is
\begin{equation*}
  g_{ij}(z,h) =  \frac{1}{ \sinh(\lambda h)|b_i|} \int_{h (a-|b_i|)}^{h (a+|b_i|)}\phi_{\sigma\sqrt{h}}(z-x) \tilde{q}_{h,-}\Big(\frac{x-ah}{b_i}\Big) dx,\quad z\in \R.
\end{equation*}
with
\begin{equation*}
  \tilde{q}_{h,-}(x) =  \frac{\lambda}2 B_0\Big(\sqrt{(h^2-x^2)} \lambda \Big), \quad x\in [-h,h].
\end{equation*}
\end{theorem}

A proof of this statement can be found in the appendix.

\subsection{Non-symmetric $Q$ matrix}
In this section we again assume that the state space consists of two elements only, i.e.  $\mathcal{E}=\{1,2\}$ but the transition intensities are arbitrary, i.e. $q_{12}=\lambda_1>0, q_{21}=\lambda_2>0$. The transition probabilities in \eqref{eq:transprob} are now given by
\begin{eqnarray}
\label{eq:12}
  p_{12}(t) &=& \frac{\lambda_1}{\lambda_1+\lambda_2}-\frac{\lambda_1}{\lambda_1+\lambda_2} e^{-(\lambda_1+\lambda_2)t}, \\
\label{eq:21}
  p_{21}(t) &=& \frac{\lambda_2}{\lambda_1+\lambda_2}-\frac{\lambda_2}{\lambda_1+\lambda_2} e^{-(\lambda_1+\lambda_2)t},
\end{eqnarray}
In this case we interpret $(I_t)$ in \eqref{telegraph} as the asymmetric telegraph process. This means that
\begin{equation*}
I_t := \int_0^t (-1)^{N_s} ds
\end{equation*}
where now
\begin{equation*}
  N_t =  |\{ 0\le s\le t : \varepsilon_s \neq \varepsilon_{s-}\}|
\end{equation*}
is the number of state transitions in $[0,t]$ of the Markov chain $(\varepsilon_t)$. Again explicit formulas for the densities $g_i, g_{ij}$, $i,j\in\{1,2\}$ can be derived. The density of $I_h$ given $\varepsilon_0=i$ has been obtained in \cite{MKR94} by using the Girsanov theorem for point processes.
In \cite{LR14} the authors use a differential equations approach. For our situation we get (w.l.o.g. we assume in the next theorem that the initial state is $\varepsilon_0=1$):

\begin{theorem}\label{theo:nonsymQ}
The density of $Z_h$ given $\varepsilon_0=1$ is
\begin{eqnarray*}
  g_1(z,h) &=&  e^{-\lambda_1 h} \Big[\phi_{\sigma\sqrt{h}}(z-h(b_1+a)) + \frac{1}{|b_1|} \int_{h (a-|b_1|)}^{h (a+|b_1|)}\phi_{\sigma\sqrt{h}}(z-x) \tilde{p}_{h}\Big(\frac{x-ah}{b_1}\Big) dx\Big],\quad z\in \R.
\end{eqnarray*}
with
\begin{eqnarray*} \tilde{p}_{h}(x) &=& \frac{\lambda_1}2 exp\Big(\frac{\lambda_1-\lambda_2}2 (h-x)\Big) \Big[B_0\Big( \sqrt{\lambda_1\lambda_2 (h^2-x^2)}\Big) \\
&&+ \Big(\frac{(h+x)\lambda_2}{(h-x)\lambda_1} \Big)^{\frac 12} B_1\Big( \sqrt{ \lambda_1\lambda_2 (h^2-x^2)}\Big)  \Big], \quad x\in [-h,h].
\end{eqnarray*}
The density of $Z_h$ given $\varepsilon_0=\varepsilon_h=1$ is for $z\in\R$
\begin{equation*}\label{eq:densityg_iisym2}
  g_{11}(z,h) = \frac{(\lambda_1+\lambda_2)}{\lambda_2 e^{\lambda_1h}+\lambda_1 e^{-\lambda_2 h}}
   \cdot \Big[ \phi_{\sigma\sqrt{h}}(z-h(b_1+a))  + \frac{1}{|b_1|} \int_{h (a-|b_1|)}^{h (a+|b_1|)}\phi_{\sigma\sqrt{h}}(z-x) \tilde{q}_{h,+,1}\Big(\frac{x-ah}{b_1}\Big) dx\Big],
\end{equation*}
with
\begin{equation*}
  \tilde{q}_{h,+,1}(x) =  \frac{\lambda_1}2 \exp\Big(\frac{\lambda_1-\lambda_2}{2}(h-x)\Big) \Big(\frac{(h+x)\lambda_2}{(h-x)\lambda_1} \Big)^{\frac 12} B_1\Big(\sqrt{\lambda_1\lambda_2(h^2-x^2)} \Big), \quad x\in [-h,h].
\end{equation*}
The density of $Z_h$ given $\varepsilon_0=1, \varepsilon_h=2$ is
\begin{equation*}
  g_{12}(z,h) =  \frac{(\lambda_1+\lambda_2)}{ e^{\lambda_1 h} - e^{-\lambda_2 h}} \cdot \frac{1}{|b_1|}\int_{h (a-|b_1|)}^{h (a+|b_1|)}\phi_{\sigma\sqrt{h}}(z-x) \tilde{q}_{h,-,1}\Big(\frac{x-ah}{b_1}\Big) dx, z\in\R.
\end{equation*}
with
\begin{eqnarray*}
  \tilde{q}_{h,-,1}(x) &=&  \frac{1}2 \exp\Big(\frac{\lambda_1-\lambda_2}{2}(h-x)\Big) B_0\Big(\sqrt{\lambda_1\lambda_2(h^2-x^2)} \Big), \quad x\in [-h,h].
\end{eqnarray*}
\end{theorem}
The proof can again be found in the appendix.

\section{Approximate filter for more states}\label{sec:approx}
Using the methods of the preceding section it is in principle possible to derive filters if $\mathcal{E}$ consists of more than two states. However, formulas get quite complicated, in particular when the state space is large. Instead we try to derive approximate filters in this case. The first approach is based on the continuous filter obtained by \eqref{eq:KS} and \eqref{eq:zakai}.

\subsection{Approach via continuous filter}
In \cite{PR10,KPS93} the solution of the Zakai equation \eqref{eq:zakai} which is a homogeneous linear It\^{o} SDE, is discussed. If the matrices $Q$ and $D$ (recall that $D$ is the diagonal matrix with elements ${\alpha_1}/{\sigma^2},\ldots {\alpha_d}/{\sigma^2}$ on the diagonal) commute, i.e. if $QD=DQ$ then an explicit solution of \eqref{eq:zakai} would be
\begin{equation}\label{eq:zakaiapprox} \xi_t = \exp\Big(Qt -\frac12 D^2 t + DZ_t \Big).\end{equation}
However this is not true if the matrices do not commute, but still the expression can serve as an approximate filter, see \cite{PR10,KPS93}. Note that this solution only involves $Z_t$ and not the path of the process and is thus also feasible for a filter with discrete observations. We will use \eqref{eq:KS} and \eqref{eq:zakaiapprox} as one approximation for our discrete filter in Section \ref{sec:numerics}.   Since this approach relies on an assumption which is not valid in general it is not possible to show any kind of convergence of this approximation scheme. Alternatively, $(\xi_t)$ can be approximated by a discrete process using the Milstein scheme in order to solve the SDE in \eqref{eq:zakai} numerically:
\begin{equation}\label{eq:Milstein}
\hat{\xi}_{n+1} = \Big[I + (Q-\frac12 D^2 ) h + D(I+\frac12 D) (Z_{t_{n+1}}-Z_{t_n}) \Big]\hat{\xi}_n.
\end{equation}
Here $\hat{\xi}_n := \hat{\xi}_{t_n}$ where $t_n = nh$ as before.
Again we see that this approximation only involves $(Z_{t_{n+1}}-Z_{t_n})$ which can be observed in our situation, too and thus the approximate filter can also be used in the discrete setting. For \eqref{eq:Milstein} it is known that the continuous approximation $\hat{\xi}$ which is obtained from $\hat{\xi}_n$ converges with strong order $1$ to the true solution of \eqref{eq:zakai} if the time step size of the discretization tends to zero. Moreover, this convergence is also true for the corresponding filter (see \cite{PR10} Theorem 4.1). Hence for frequent observations the approximation in \eqref{eq:Milstein} is supposed to work well in our setting. Note however, that a practical implementation of \eqref{eq:Milstein} requires some smoothing in order to avoid affects like $\hat{\xi}_n$ getting negative.

\subsection{Approach via PDEs}
Here we derive a system of PDEs for the densities ${g}_{ij}$  which appear in the filter \eqref{filter}. From $g_{ij}$ we obtain the densities $g_i$ by
\eqref{eq:gifromgij}. As soon as continuous-time Markov chains are involved it is a common tool to work with differential equations, see e.g. the system of ODEs for the probability distribution of a Markov chain \eqref{eq:Q} or the partial differential equation for the position of a particle given by the telegraph process, see \cite{Kac}. This equation is a hyperbolic second order differential equation and is known as the telegraph (or damped wave) equation. Partial differential equations for the densities of the telegraph process in the asymmetric case have been derived in \cite{LR14}.

In what follows define for $i,j\in\mathcal{E}$, $x\in\R$ and $t\ge 0$ the distribution
\begin{equation}
F_{ij}(x,t) = \Pop_i(J_t \le x, \varepsilon_t=j)
\end{equation}
and the corresponding density by $\tilde{f}_{ij}(x,t)$. By $f_{ij}(x,t)$ we denote the density of $\Pop(J_t \le x |\varepsilon_0=i, \varepsilon_t=j)$.
Thus, we obviously have the following relation:
\begin{equation}\label{eq:fij_from_pde}
\tilde{f}_{ij}(x,t) = {f}_{ij}(x,t) p_{ij}(t).
\end{equation}
From $f_{ij}$ we can compute immediately by convolution the required density $g_{ij}$ for our filter using \eqref{eq:Z}. Note that the density $\tilde{f}_{ii}(x,t)$ has an atom at $x=\alpha_i t$ of size $e^{-q_i t}$ where $q_i = -q_{ii}$, since with probability $e^{-q_i t}$ no switch occurs in the time interval $[0,t]$ in which case $J_t =\alpha_i t$.

We will now derive a PDE for $\tilde{f}_{ij}$ by conditioning on what happens in the first time interval $[0,h]$. We obtain
\begin{eqnarray*}
F_{ij}(x,t) &=&  (1-q_i h +o(h)) F_{ij}(x-\alpha_i h, t-h) +\\
  && \int_0^h \sum_{k\neq j} q_i e^{-q_i s} p_{ik} F_{ik}(x-\alpha_i s,t-s) ds + o(h),
\end{eqnarray*}
where $p_{ik}$ is the transition probability of the embedded Markov chain given by
$$ p_{ik} = \frac{q_{ik}}{q_i},\quad\mbox{for}\; i\neq k.$$
Rearranging terms and dividing by $h$ yields:
\begin{eqnarray*}
&& \frac1h (F_{ij}(x,t)-F_{ij}(x-\alpha_ih,t)) + \frac1h (F_{ij}(x-\alpha_ih,t)-F_{ij}(x-\alpha_i h,t-h)) =\\
 &=&  -q_i F_{ij}(x, t) +\frac1h  \int_0^h \sum_{k\neq i} q_i e^{-q_i s} p_{ik} F_{ik}(x-\alpha_i s,t-s) ds + \frac{o(h)}h.
 \end{eqnarray*}
Now letting $h\downarrow 0$ we obtain (note that the limit on the right hand side exists, hence also on the left hand side):
\begin{eqnarray*}
&& \alpha_i \frac{\partial}{\partial x}F_{ij}(x,t) +  \frac{\partial}{\partial t}F_{ij}(x,t)  = \sum_k q_{ik} F_{kj}(x,t).
\end{eqnarray*}
Applying  the operator $\frac{\partial}{\partial x}$ on both sides, we derive the following equation for the densities:
\begin{eqnarray*}
&& \alpha_i \frac{\partial}{\partial x}\tilde{f}_{ij}(x,t) +  \frac{\partial}{\partial t}\tilde{f}_{ij}(x,t)  = \sum_k q_{ik} \tilde{f}_{kj}(x,t).
\end{eqnarray*}
In matrix form this equations can be written as
\begin{equation}\label{eq:pdeSystemMatrixForm}
D^\alpha \frac{\partial}{\partial x} G(x,t) + \frac{\partial}{\partial t} G(x,t) = Q G(x,t)
\end{equation}
where $G=(\tilde{f}_{ij})_{i,j\in\mathcal{E}}$ and $D^\alpha$ is a diagonal matrix with $\alpha_1,\ldots,\alpha_n$ on the diagonal. Moreover, we have the boundary conditions $\tilde{f}_{ij}(x,0) = \delta_{ij}$. This system of PDEs can then be solved numerically. The domain of $\tilde{f}_{ij}(x,t)$ is given by $x\in[\min_i \alpha_i t, \max_i \alpha_i t]$ and $t>0$. Outside the densities vanish.
This approach is exact, however involves a numerical computation. One difficulty is the point mass which appears. However, the point mass is known and can be separated. More precisely, we know that $\tilde{f}_{ii}(x,t)= \delta_{\alpha_i t}(x) e^{-q_i t} + \hat{f}_{ii}(t,x)$ where $\hat{f}$ is smooth.

\subsection{Approach via Discretization}
Here we use a very simple approximation of the densities involved in the filter by discretizing the continuous-time Markov chain. More precisely, when we consider $(\varepsilon_{m\cdot \Delta t})_m$ with $\Delta t = h/N$, then we obtain a discrete-time Markov chain with transition matrix $p(\Delta t) = e^{Q \Delta t}$. Note that for $i\neq j$ either $p_{ij}(t)\equiv 0$ or $p_{ij}(t)>0$ for all $t>0$. We assume that the latter case is valid for all $i\neq j$.  We can now approximate the random variable $J_h$ in \eqref{eq:Jt} by
\begin{equation*}
\hat{J}_h := \Delta t \sum_{m=1}^N  \alpha_{\varepsilon_{m\cdot \Delta t}}.
\end{equation*}
Since $\mathcal{E}$ is finite, the random variable $\hat{J}_h$ is obviously discrete and can only take a finite number of values. We denote by $\mathcal{D}$ the finite set of all possible values. Furthermore, we denote a multi-index $\mathbf{j}$ as a vector $\mathbf{j}=(j_1,\ldots,j_N)$, where each $j_k \in \mathcal{E}$, i.e. $\mathbf{j}\in \mathcal{E}^N$. For $d\in \mathcal{D}$ denote
\begin{equation*}
\Lambda(d) := \{ \mathbf{j}\in \mathcal{E}^N : \alpha_{j_1}+\ldots + \alpha_{j_N}= d/\Delta t\}.
\end{equation*}
Then we obtain
\begin{eqnarray*}
\nonumber \Pop_i( \hat{J}_h = d) &=& \Pop_i\Big(  \sum_{m=1}^N  \alpha_{\varepsilon_{m\cdot \Delta t}} = d/\Delta t\Big) \\
   &=& \sum_{\mathbf{j}\in \Lambda(d)} \Pop_i(\varepsilon_{\Delta t}=j_1) \prod_{k=1}^{N-1} \Pop( \varepsilon_{(k+1)\Delta t} = j_{k+1} | \varepsilon_{k \Delta t} = j_k).
\end{eqnarray*}
Thus, we can approximate the density $g_i(z,h)$ of the distribution of $Z_h$ given $\varepsilon_0=i$ by
\begin{eqnarray}\label{eq:giGMApprox}
 \hat{g}_i(z,h) &=& \sum_{d\in \mathcal{D}} \Pop_i(  \hat{J}_h = d) \phi_{\sigma\sqrt{h}} (z-d).
\end{eqnarray}
Analogously, for fixed $i,j\in \mathcal{E}$ and $d\in \mathcal{D}$ denote
\begin{equation*}
\Lambda(j,d) := \{ \mathbf{j}\in \mathcal{E}^{N-1} : \alpha_{j_1}+\ldots + \alpha_{j_{N-1}}= d/\Delta t - \alpha_j\}.
\end{equation*}
Then we obtain
\begin{eqnarray*}
\nonumber &&\Pop_{i}( \hat{J}_h = d|\varepsilon_h=j) = \Pop_{i}\Big(  \sum_{m=1}^{N-1}  \alpha_{\varepsilon_{m\cdot \Delta t}} = d/\Delta t-\alpha_j\Big| \varepsilon_h=j\Big) \\
   &=& \sum_{\mathbf{j}\in \Lambda(j,d)} \frac{\Pop_i(\varepsilon_{\Delta t}=j_1)}{\Pop_i(\varepsilon_h = j)} \prod_{k=1}^{N-2} \Pop( \varepsilon_{(k+1)\Delta t} = j_{k+1} | \varepsilon_{k \Delta t} = j_k) \cdot \Pop( \varepsilon_{h} = j | \varepsilon_{(N-1) \Delta t} = j_{N-1}) .
\end{eqnarray*}
Thus, we can approximate the density $g_{ij}(z,h)$ of the distribution of $Z_h$ given $\varepsilon_0=i, \varepsilon_h=j$ by
\begin{eqnarray}\label{eq:gijGMApprox}
 \hat{g}_{ij}(z,h) &=& \sum_{d\in \mathcal{D}} \Pop_{i,j}(  \hat{J}_h = d) \phi_{\sigma\sqrt{h}} (z-d).
\end{eqnarray}
It is rather obvious that in this approach the approximate density converges against the exact one if $\Delta t\to 0$.

\section{Numerical Examples}\label{sec:numerics}
This section illustrates the results presented in this work using some numerical examples. Basically, we consider two cases. First, we evaluate a scenario in which the underlying continuous time Markov chain has two states in order to provide an illustration comparing the exact filter to the approximations.  Second, we consider a five state case. For application of the continuous filter \eqref{eq:zakaiapprox}, it is important to keep in mind, that the induced approximation error not only stems from the fact that the we have discrete measurements, but also from the fact that the considered matrices $D$, and $Q$ do not commute.

All examples were implemented using Matlab R2014a on a laptop with an Intel i7-2620M CPU and 8GB RAM. The system of PDEs \eqref{eq:pdeSystemMatrixForm} was solved using the algorithm proposed by \cite{Skeel1990}, which is directly implemented in Matlab. From the solution $\tilde{f}_{ij}$, we derived $f_{ij}(x,t)$ using \eqref{eq:fij_from_pde} 
Then, a numerical convolution procedure was performed for obtaining the values of the density $g_{ij}(x,t)$ on the grid points and $g_{i}(x,t)$ was obtained using \eqref{eq:gifromgij}. Cubic interpolation was used for evaluating these densities at other points.

\subsection{Two States Example}
\begin{figure}
\includegraphics[width=0.48\textwidth]{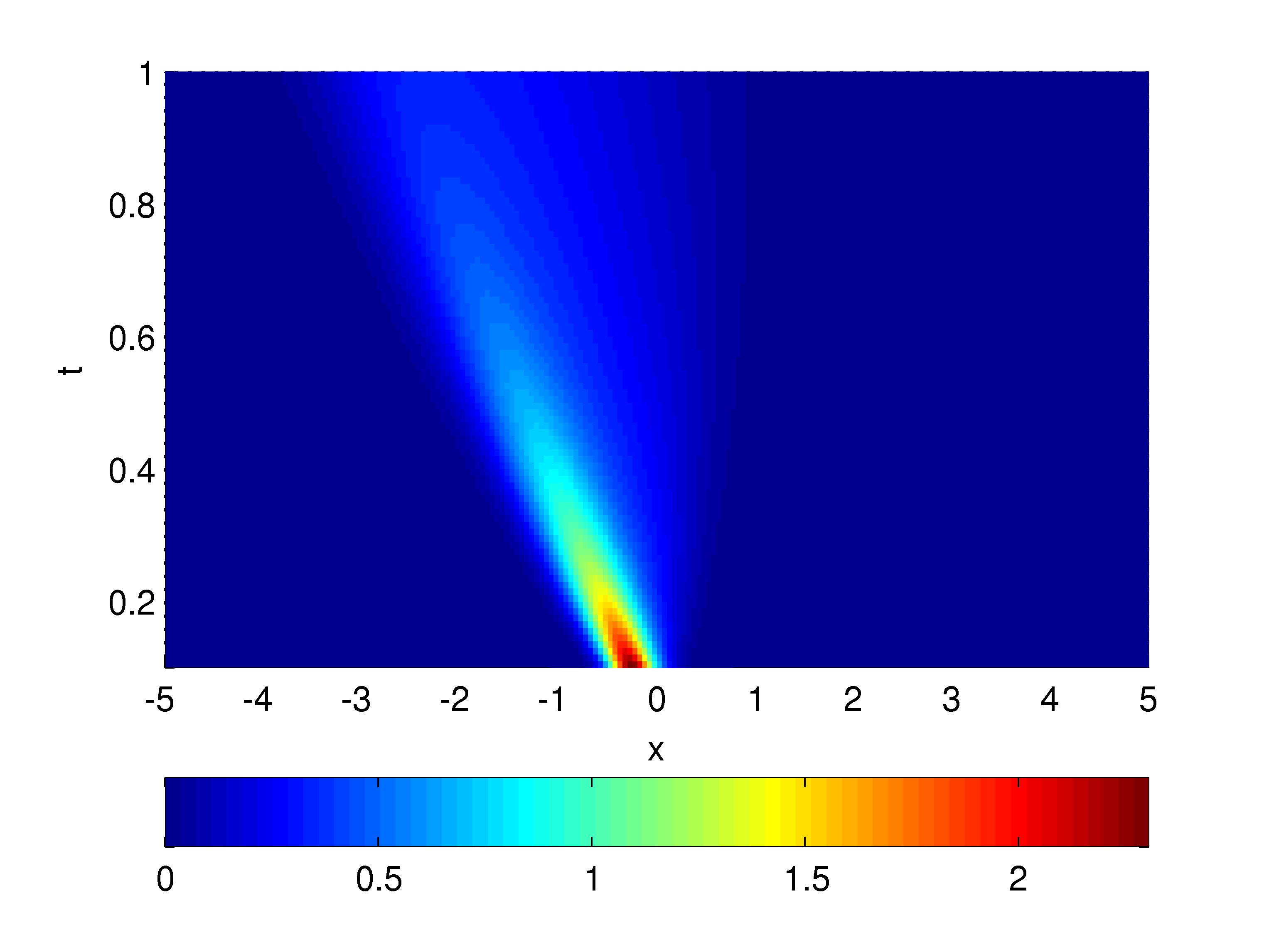}
\hfill
\includegraphics[width=0.48\textwidth]{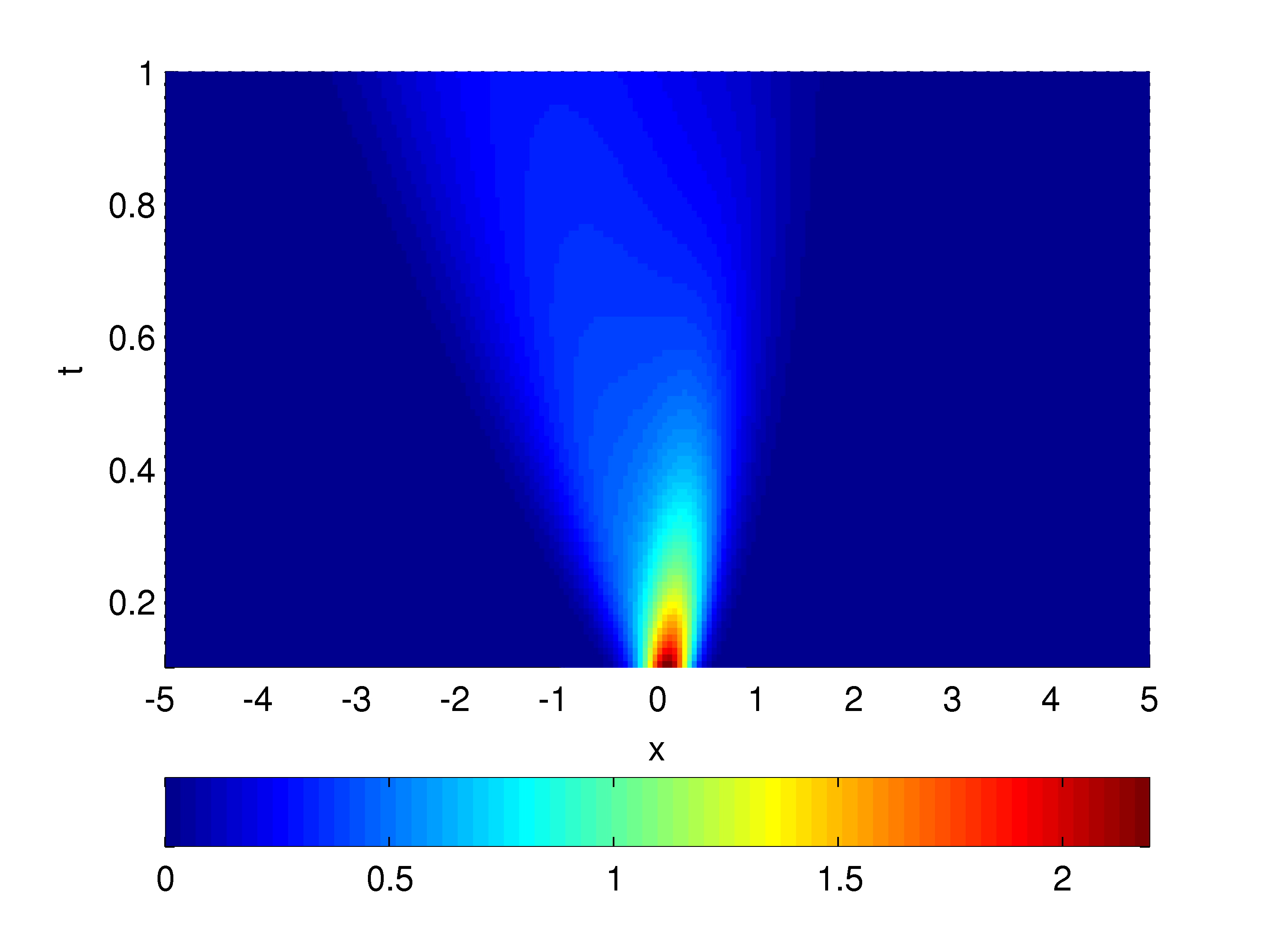}
\caption{Density of $Z_t$ conditioned on $\epsilon_0=1$ (left) and $\epsilon_0=2$ (right).\label{fig:densityStartCond}}
\end{figure}

\begin{figure}
\subfigure[$\varepsilon_0=1,\,\varepsilon_t=1$]{\includegraphics[width=0.48\textwidth]{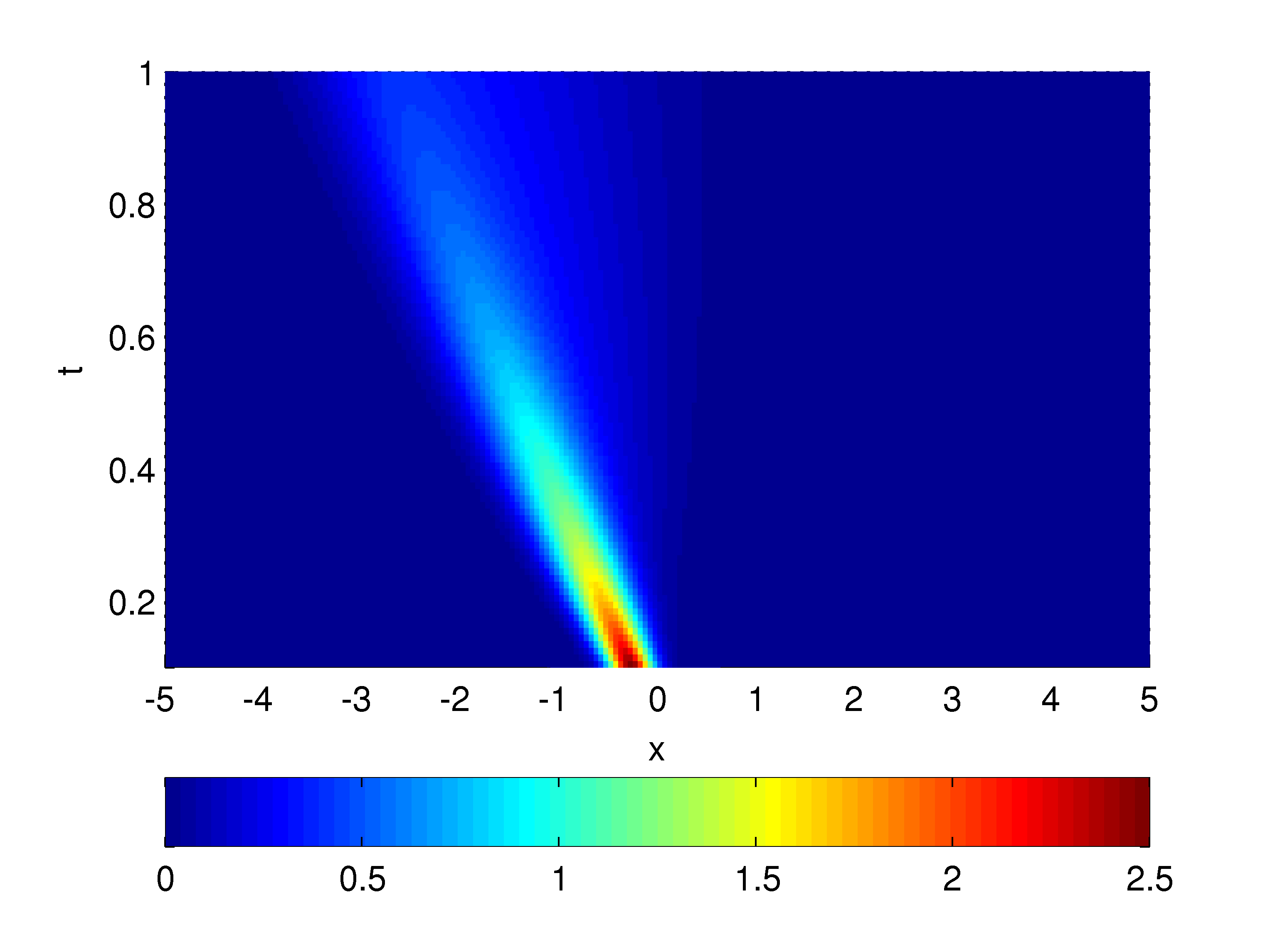}}
\hfill
\subfigure[$\varepsilon_0=1,\,\varepsilon_t=2$]{\includegraphics[width=0.48\textwidth]{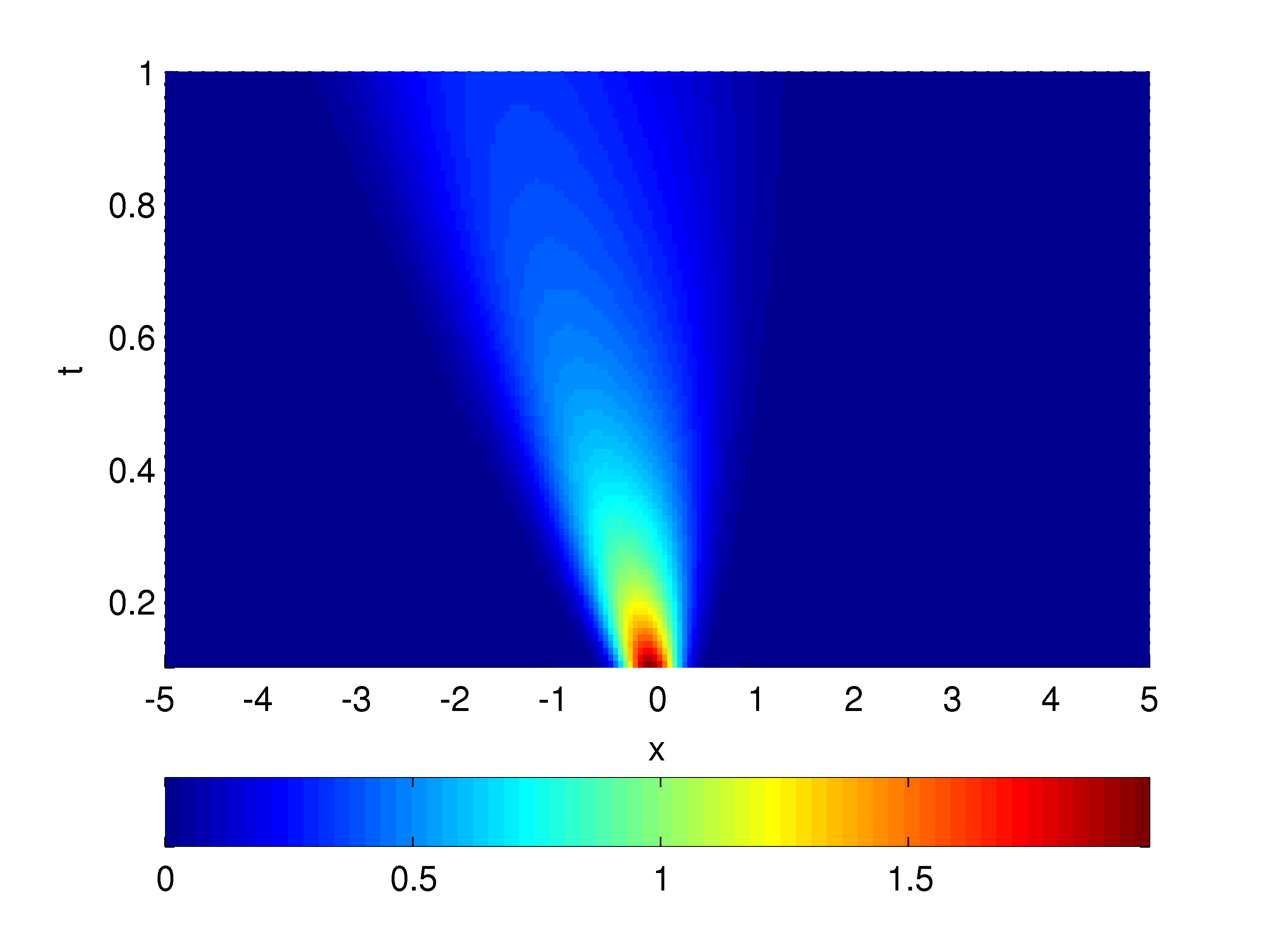}}\\
\subfigure[$\varepsilon_0=2,\,\varepsilon_t=1$]{\includegraphics[width=0.48\textwidth]{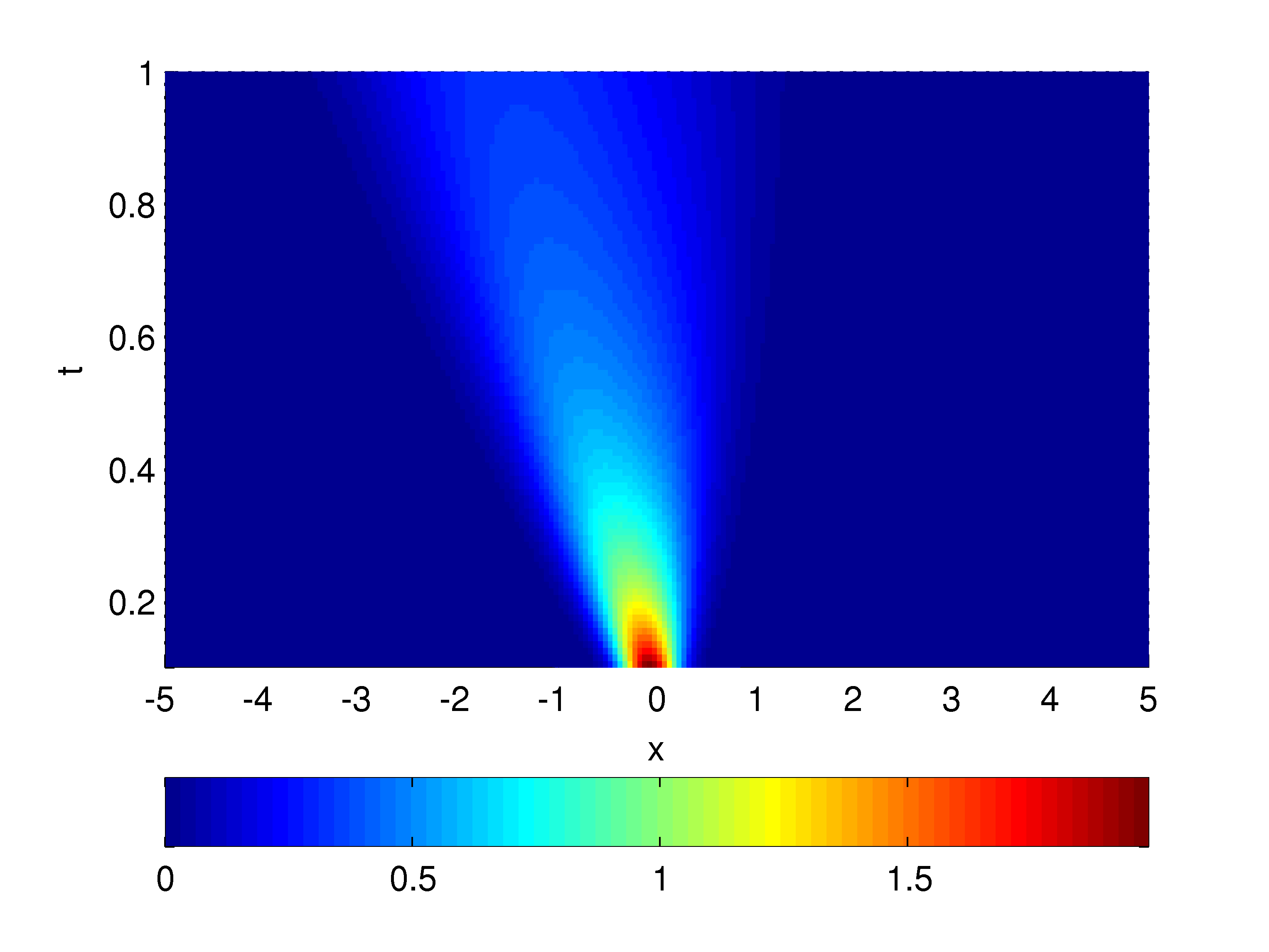}}
\hfill
\subfigure[$\varepsilon_0=2,\,\varepsilon_t=2$]{\includegraphics[width=0.48\textwidth]{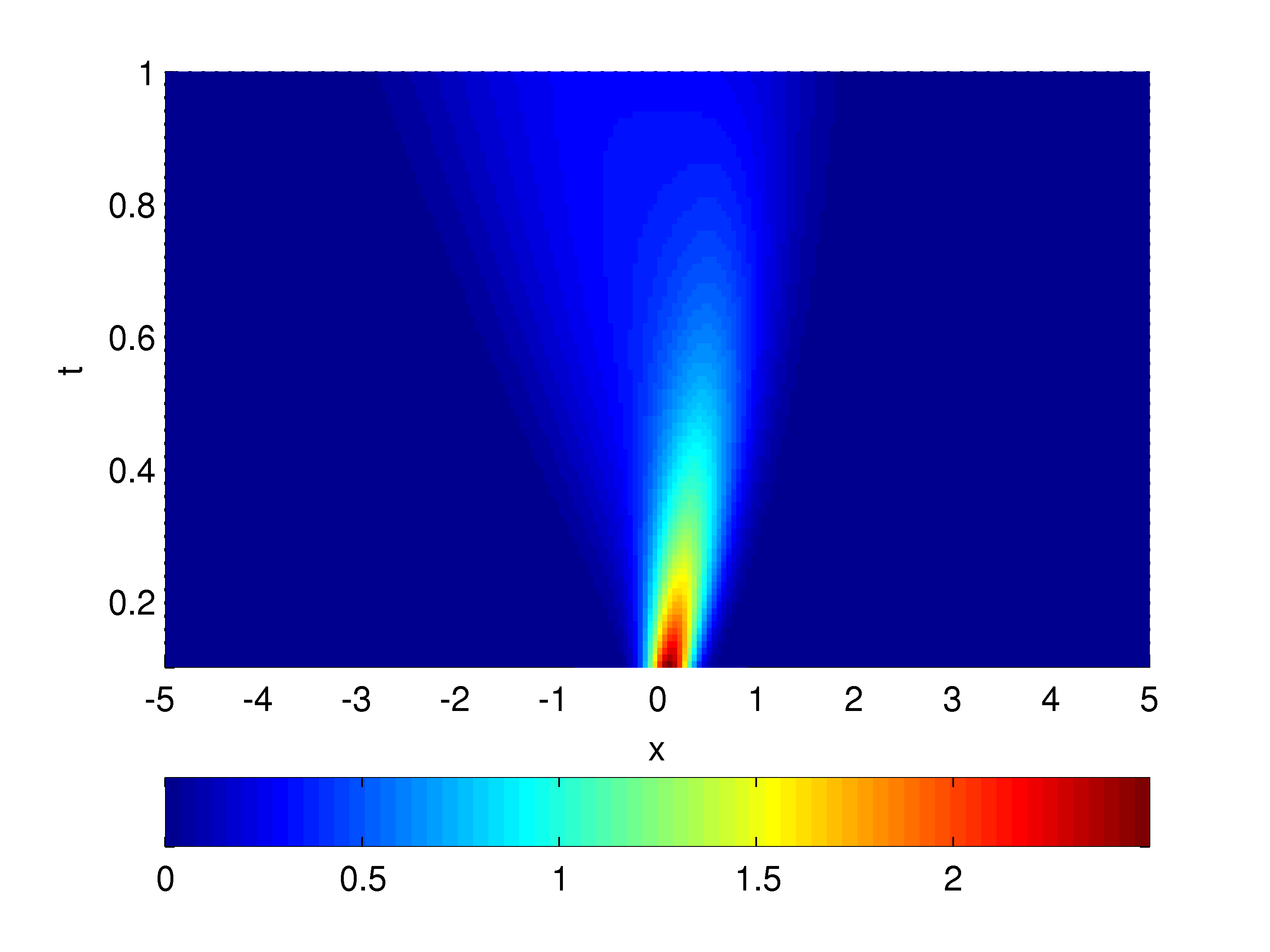}}
\caption{Density of $Z_t$ conditioned on $\varepsilon_0$ and $\varepsilon_t$.\label{fig:densityBothCond}}
\end{figure}

\begin{figure}\centering
\includegraphics[width=0.8\textwidth]{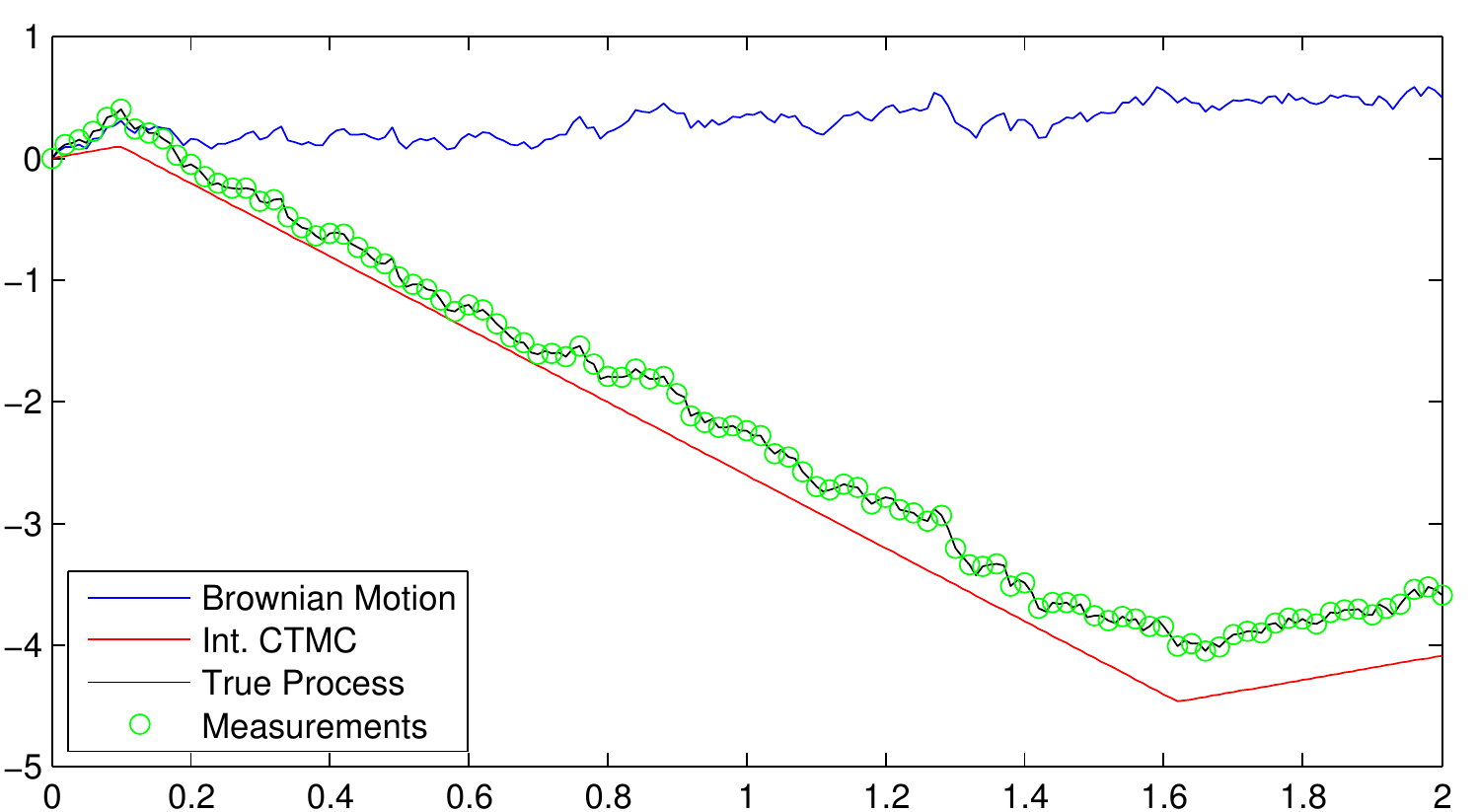}
\caption{Observed process and its components (2-state case).\label{fig:twoStateGT}}
\end{figure}

We consider the scenario $\alpha_1=-3$, $\alpha_2=1$, $\sigma=1$, initial distribution of the continuous-time Markov chain $p_0=(0.1,\ 0.9)$, and intensity matrix
\begin{equation*}
Q = \begin{pmatrix}
-2 & 2\\
 3 &-3
\end{pmatrix}.
\end{equation*}
In fig.~\ref{fig:densityStartCond}, we show the density of $Z_t$ conditioned on the value of $\epsilon_0$. In fig.~\ref{fig:densityBothCond}, we show the density of $Z_t$ conditioned on both, $\epsilon_0$ and $\epsilon_t$. It can be seen that the additional restriction of $\epsilon_t$ has quite some impact on the shape of the density. For evaluation of the filters, the observed process was generated using 100 discretization points and 50 observations in each time step. Both, the process and its observations are shown in fig.~\ref{fig:twoStateGT}

\begin{table}\centering
\subtable[Two States Example]{
\begin{tabular}{|l|c|}
\hline
{\bf Filter} &  {\bf Time}\\ \hline
Exact & 248.53 \\ \hline
PDE Based & 3.32 \\ \hline
Discretized & 4.73\\ \hline
Quasi-Exact & 0.18\\ \hline
\end{tabular}
}\hspace{40mm}
\subtable[Five States Example]{
\begin{tabular}{|l|c|}
\hline
{\bf Filter} &  {\bf Time}\\ \hline
Exact & - \\ \hline
PDE Based & 4.99 \\ \hline
Discretized & 171.48\\ \hline
Quasi-Exact & 0.23\\ \hline
\end{tabular}
}
\caption{Computation Time (in ms) per Filter Step}
\label{tab:RunTimes}
\end{table}

\begin{figure}\centering
\subfigure[Every measurement used]{\includegraphics[width=0.8\textwidth]{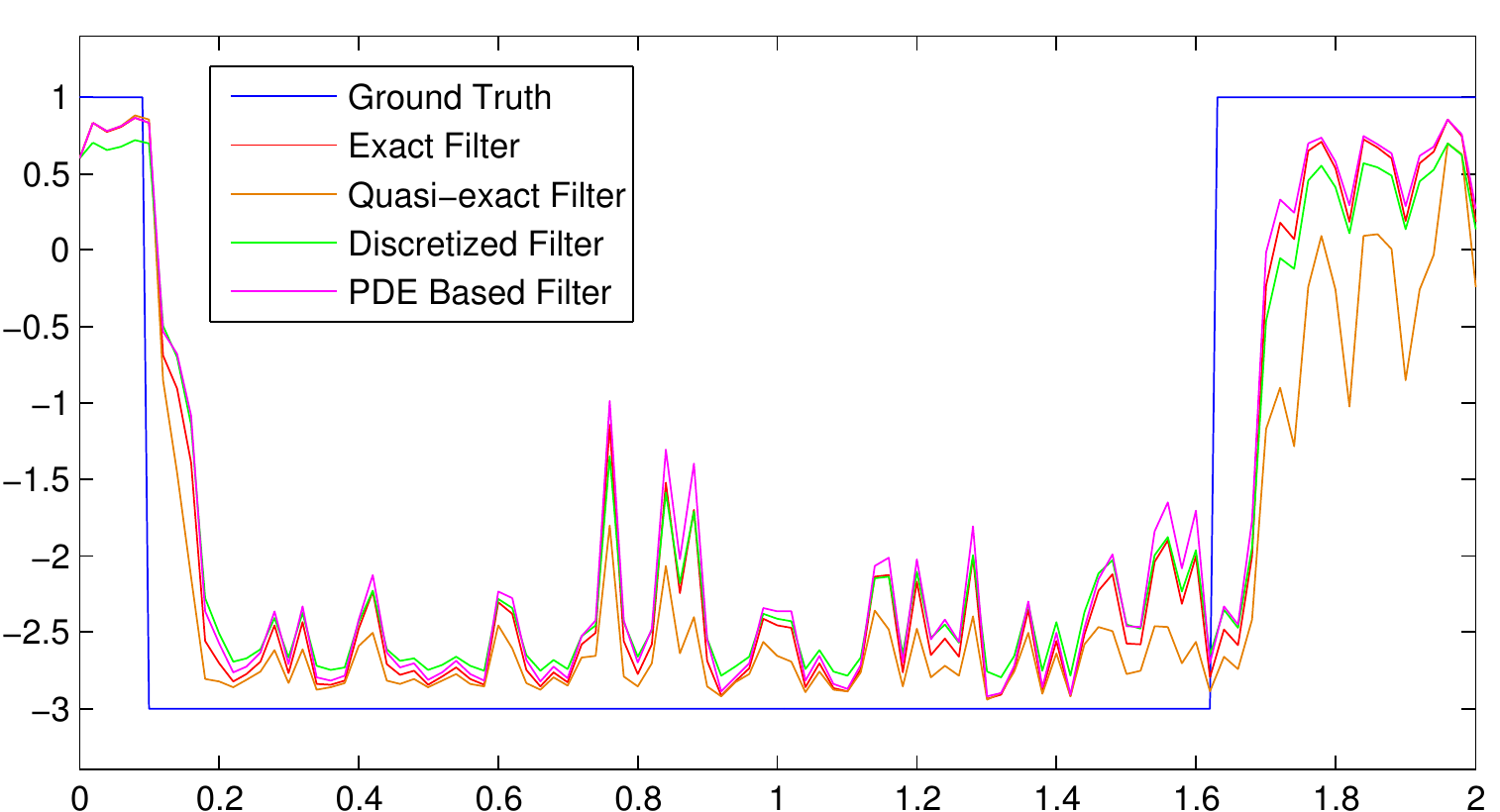}}
\subfigure[Every fifth measurement used]{\includegraphics[width=0.8\textwidth]{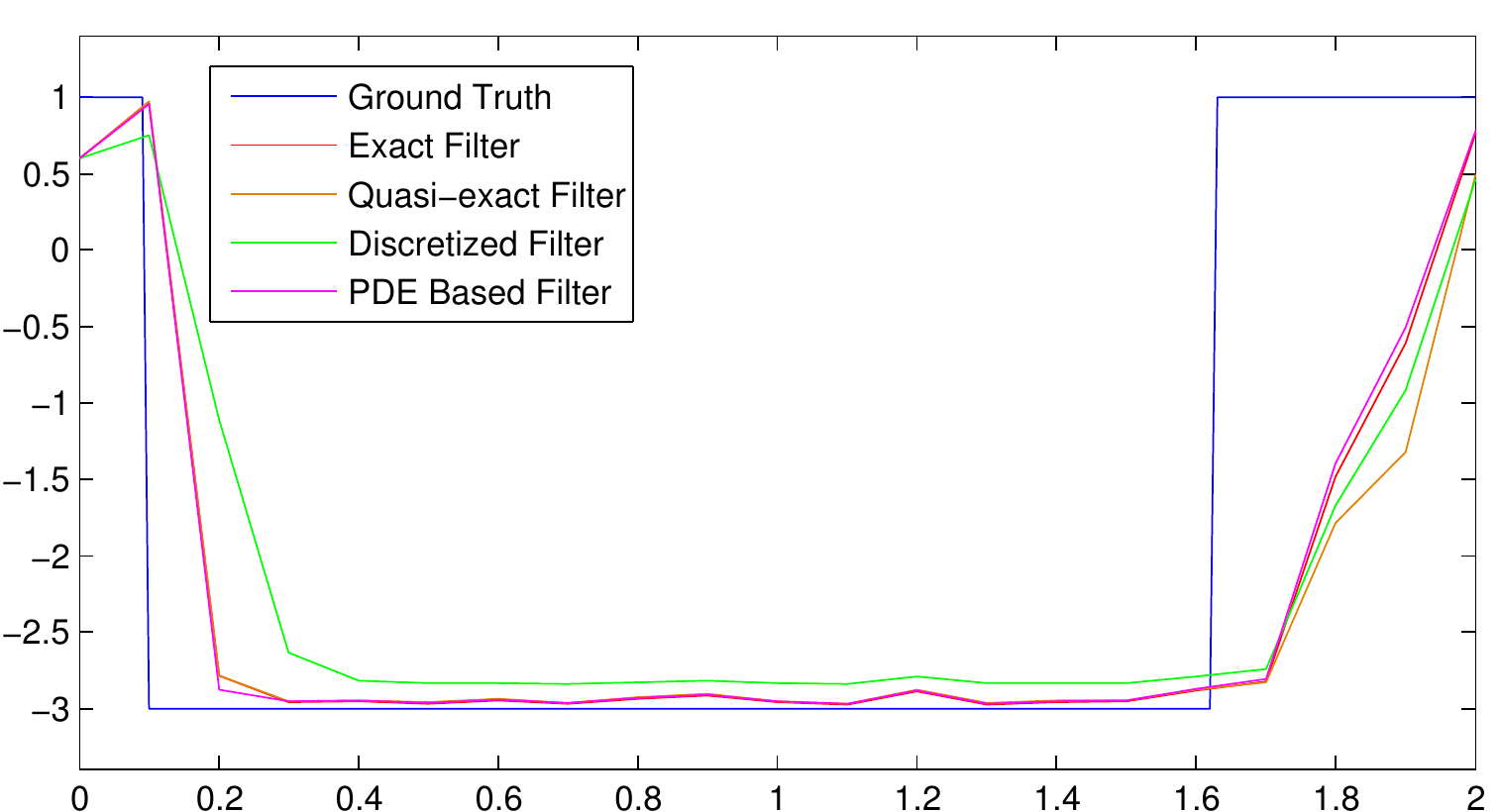}}
\caption{True state of CTMC and its estimates (2-state case).\label{fig:twoStateFilter}}
\end{figure}

In fig.~\ref{fig:twoStateFilter}, we show filter results in terms of the expectation value of the true state. This is done for two scenarios. In the first case, every measurement is used. In the second case, only every fifth measurement is used. As expected, more frequent measurements result in a faster convergence of the estimate after a change of the underlying state. On the other hand, the filter with fewer observations is less sensitive to observed errors. The exact filter is almost indistinguishable from the PDE based approach. The time discretization based filter assumed the Markov process to perform at most one jump between each measurement, i.e., we used $N=1$.

The computation time of the filters are given in table \ref{tab:RunTimes}. The high computation time for the exact filter is due to the need for performing numerical integration in each filter step. The numbers of the PDE based filter do not include the time required for numerically solving the PDE, because the resulting solution is reused in every filter step, and thus, the impact on the average computation time depends on the number of performed filter steps. The PDE was solved on an equidistant grid on $[-5,5]$ consisting of 3000 discretization points for 2000 equidistantly distributed time-points on $[0,1]$.

\subsection{Five States Example } 
\begin{figure}\begin{center}
\subfigure[$\sigma=1$]{\includegraphics[width=0.8\textwidth]{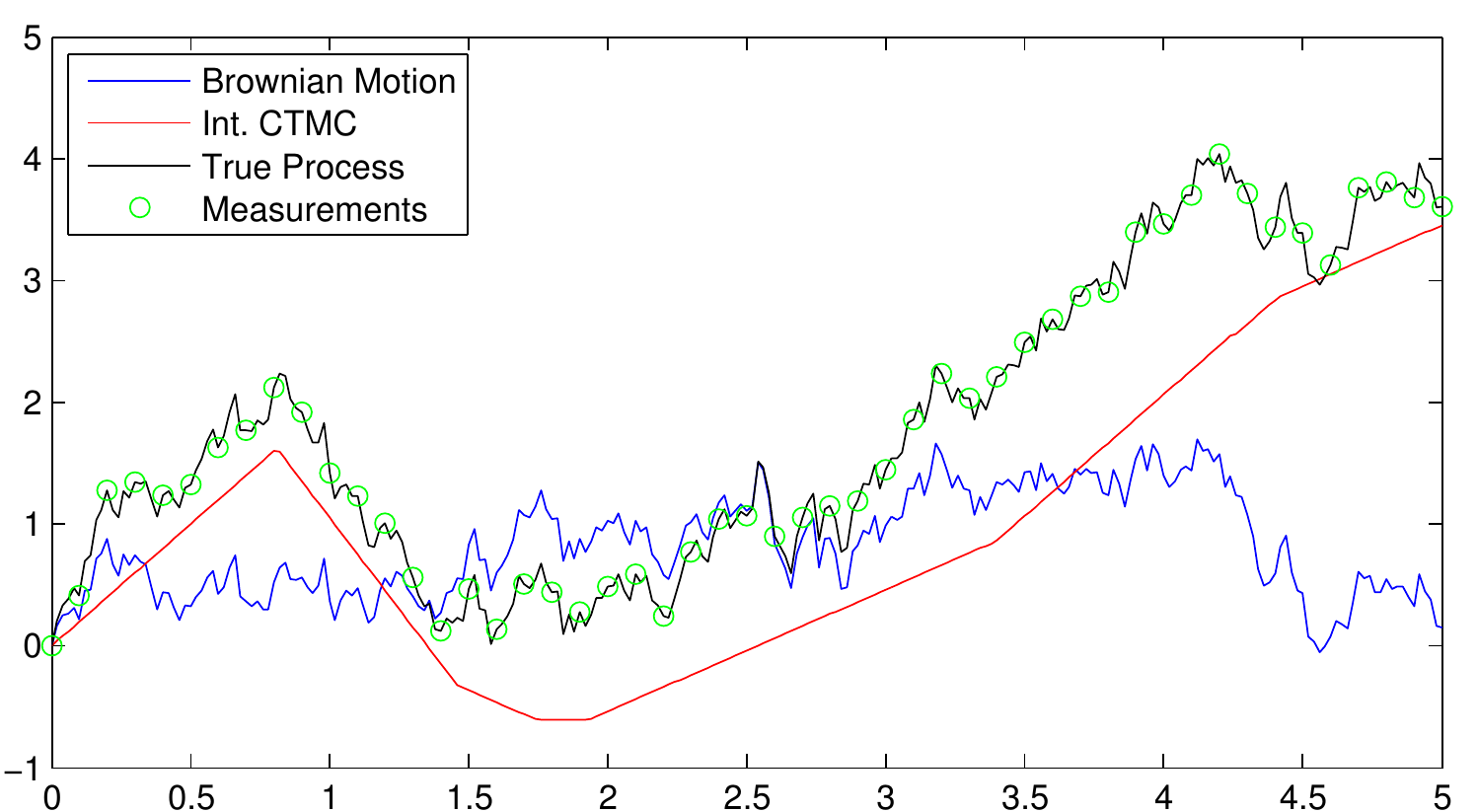}}
\subfigure[$\sigma=2$]{\includegraphics[width=0.8\textwidth]{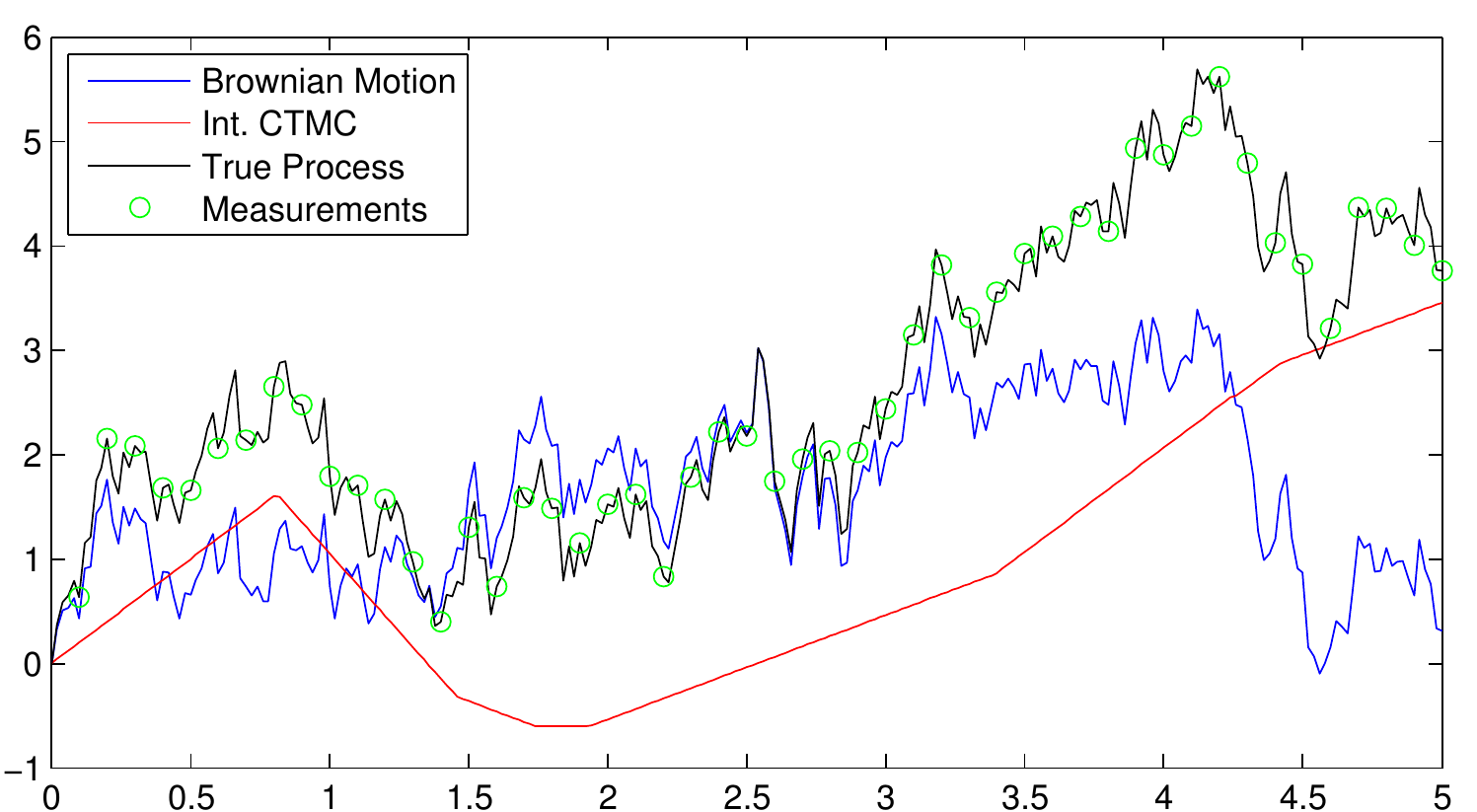}}
\caption{Observed process and its components (5-state case).\label{fig:multiStateGT}}
\end{center}\end{figure}

Now, we take a look at a more complex example where the $\alpha_i$ (for $i=1,...,5$) are given by the values $-3$, $-1$, $0$, $1$, $2$ respectively. Furthermore, we used $p_0=(0.1,\ 0.3,\ 0.3,\ 0.2,\ 0.1)$ and
\begin{equation*}
Q =
\begin{pmatrix}
 -1 & 0.5 & 0.3 & 0.1 & 0.1\\
0.4 &  -1 & 0.3 & 0.1 & 0.2\\
0.1 & 0.1 &  -1 & 0.4 & 0.4\\
0.1 & 0.1 & 0.3 &  -1 & 0.5\\
0.1 & 0.1 & 0.3 & 0.5 &  -1
\end{pmatrix}
\end{equation*}
as initial distribution and intensity matrix. Two different diffusion parameters, $\sigma=1$ and $\sigma=2$ were considered in this example in order to observe the impact of $\sigma$ on the quality of the filters. This time, we simulated 5 time steps with 50 discretization points and 10 observations per time step. The true process is shown in fig.~\ref{fig:multiStateGT}.

\begin{figure}
\subfigure[$\sigma=1$, every measurement used]{\includegraphics[width=0.48\textwidth]{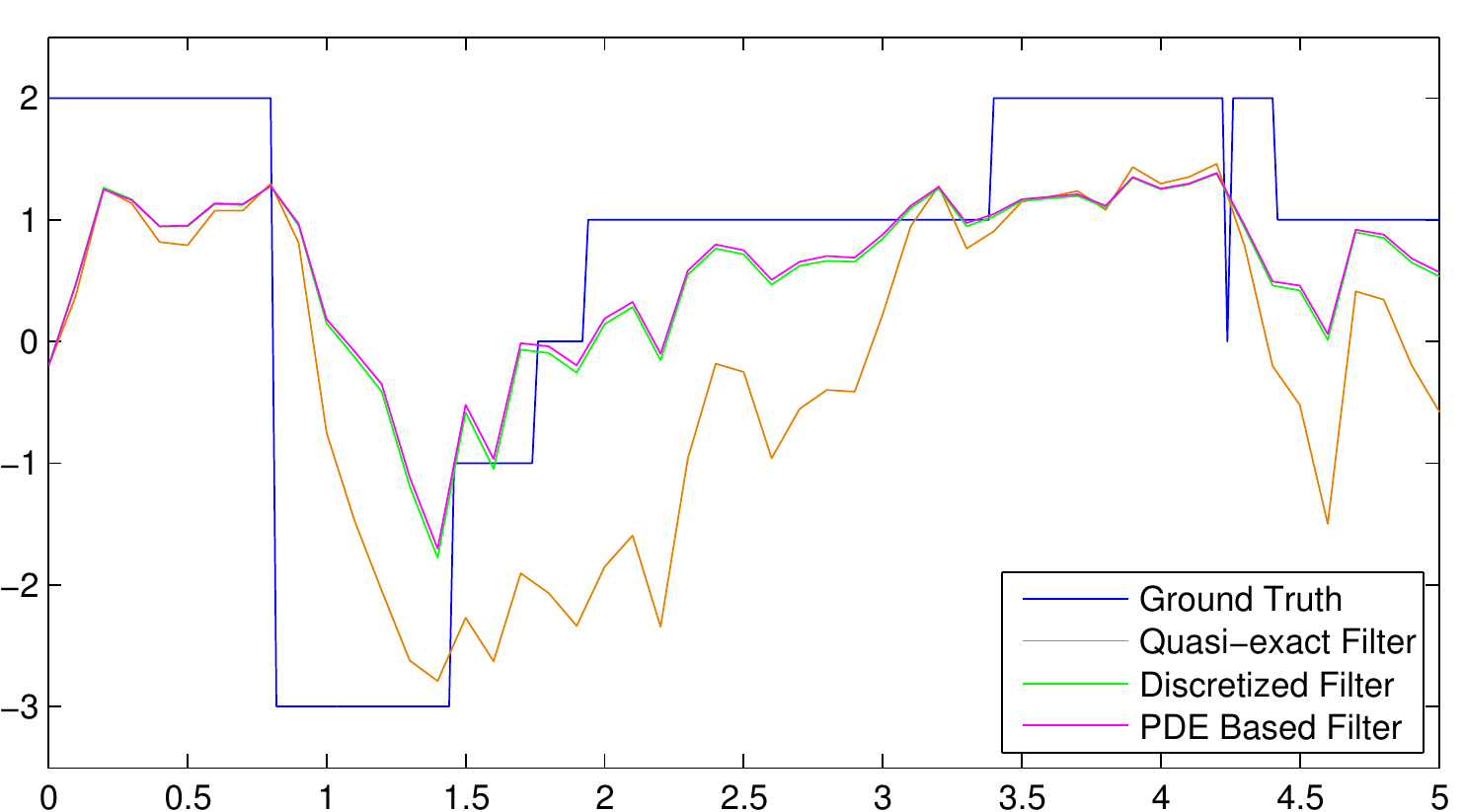}}\hfill
\subfigure[$\sigma=2$, every measurement used]{\includegraphics[width=0.48\textwidth]{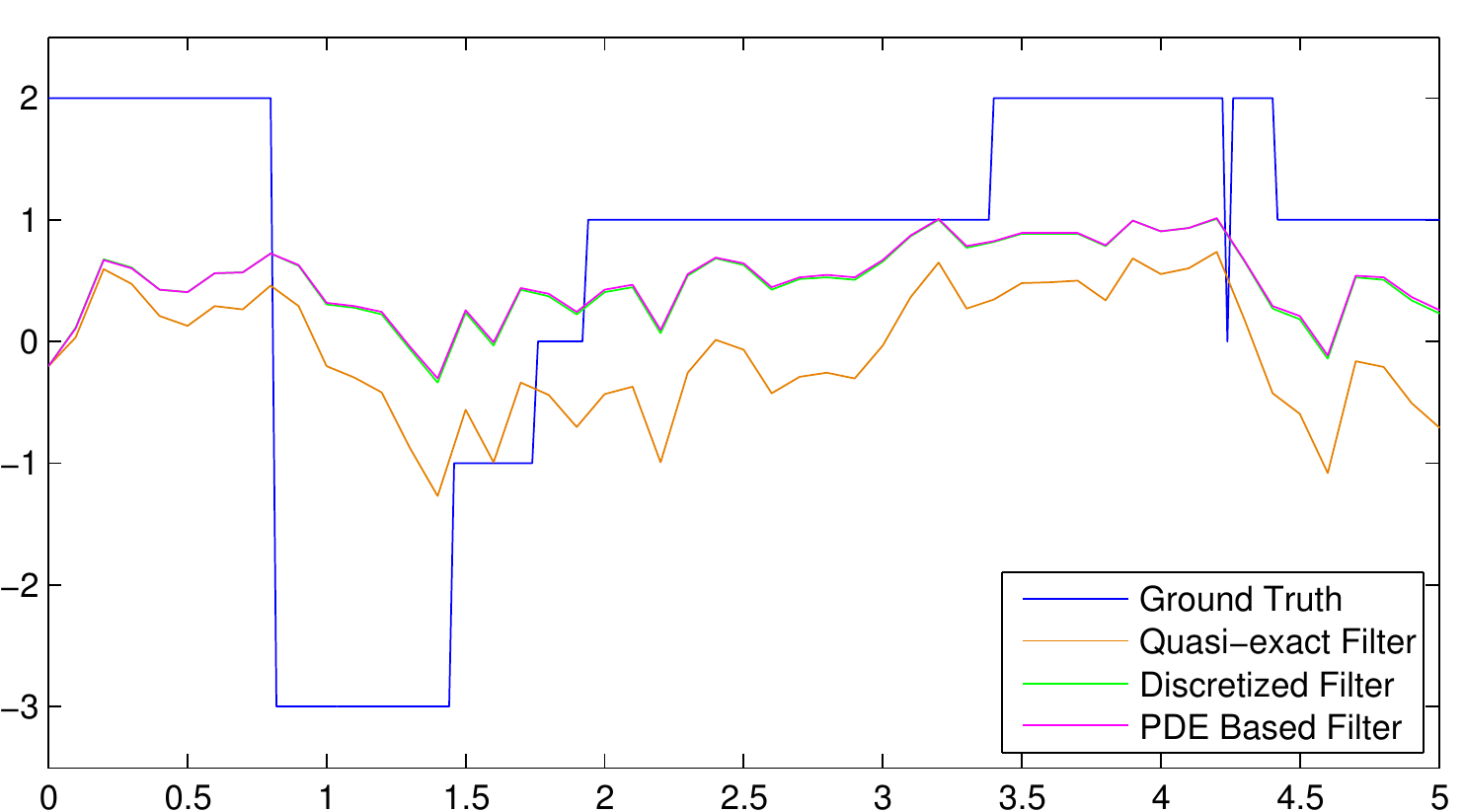}}
\subfigure[$\sigma=1$, every fifth measurement used]{\includegraphics[width=0.48\textwidth]{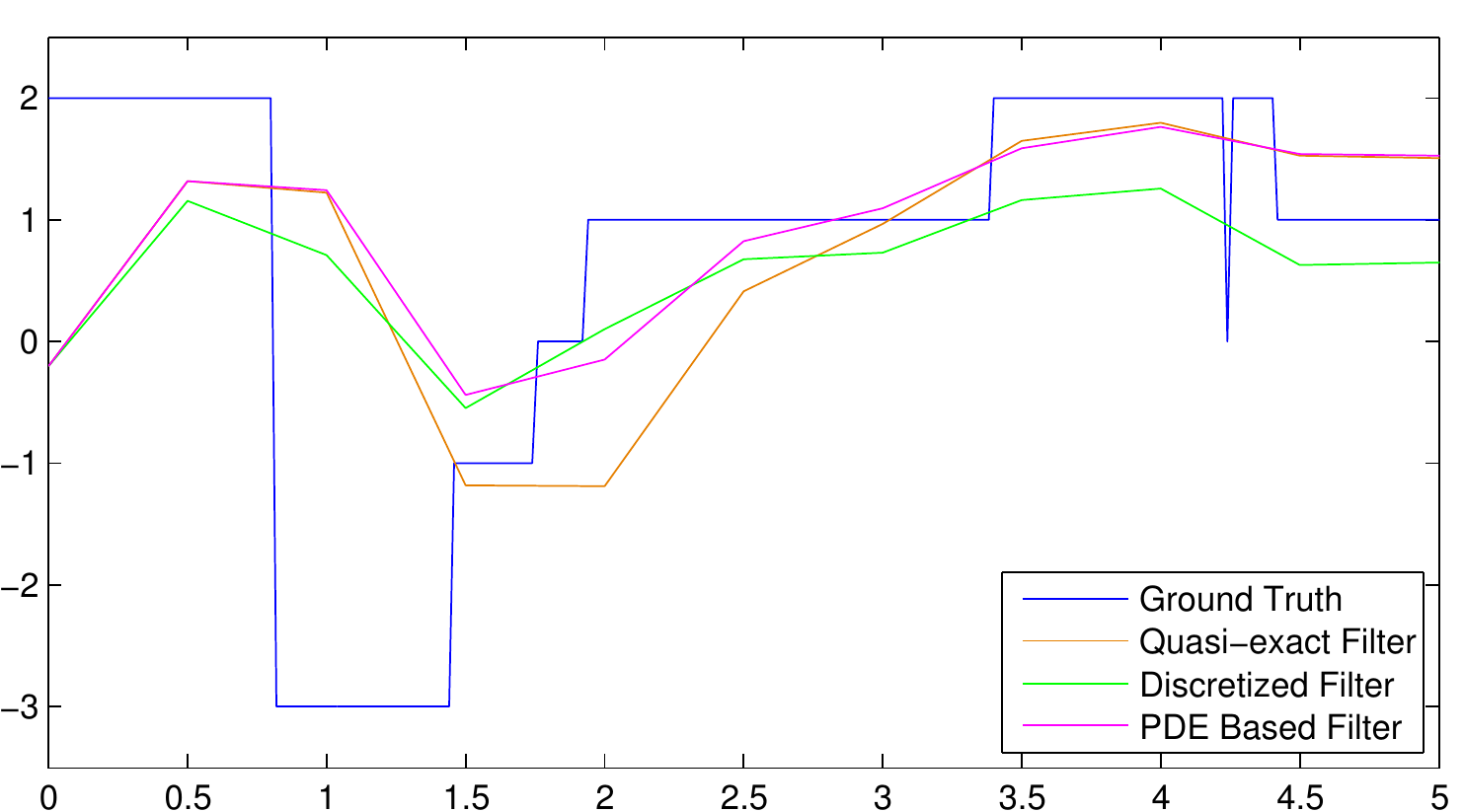}}\hfill
\subfigure[$\sigma=2$, every fifth measurement used]{\includegraphics[width=0.48\textwidth]{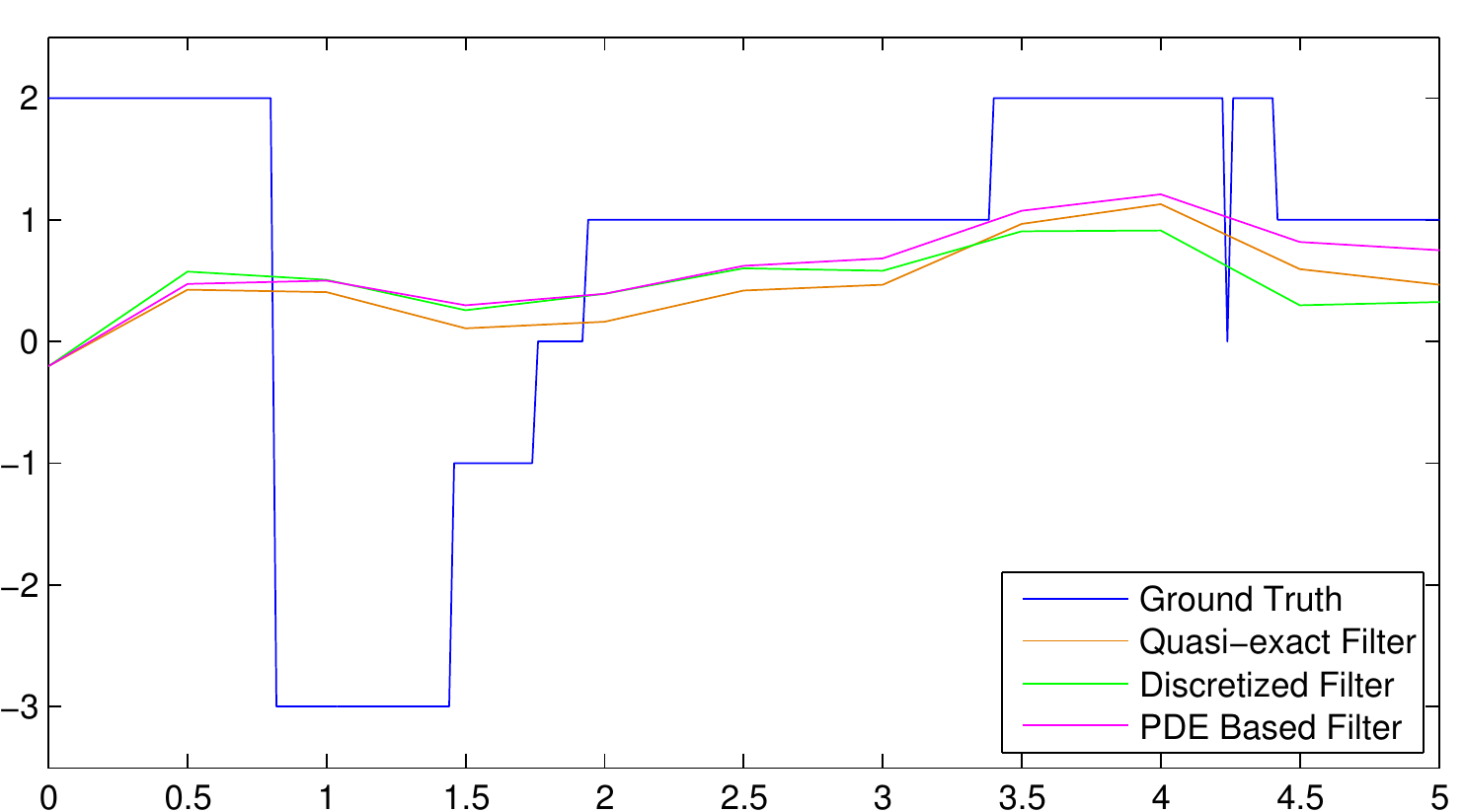}}
\caption{True state of CTMC and its estimates (5-state case).\label{fig:multiStateFilter}}
\end{figure}

The filter results are shown in fig.~\ref{fig:multiStateFilter} (once again by showing the expectation value of the obtained distributions). For a high number of measurements, the PDE based filter is almost indistinguishable from the discretized filter (which uses 4 discretization cells here, i.e., $N=4$). This is due to the fact that the PDE based approach yields the exact density which is only corrupted by  numerical errors induced from solving the PDE and by interpolation errors arising when evaluating the densities. Obviously, the expected number of jumps between two consecutive measurements decreases as the frequency of measurements goes up. Thus, in the considered scenario, the densities obtained from the discretization approach yield a good approximation of the exact density when every measurement is used. This approximation quality does not depend on the diffusion parameter $\sigma$. However, a larger $\sigma$ results in less clear differences between all filtering approaches, because less information can be obtained about the underlying true state of the continuous time Markov chain.

Looking at the computation time in table \ref{tab:RunTimes} gives a somewhat different picture this time compared to the two states case. Once again, we only show the computation times of the run using every measurement and $\sigma=1$ (it does not differ significantly for other runs). This time, the discretized filter requires a significantly higher amount of computation time. This is due to its recursive implementation and use of a higher number of discretization cells. However, this could also be optimized by obtaining the set $\mathcal{D}$ and the densities $\hat{g}_{i}(z,h)$ and $\hat{g}_{ij}(z,h)$ (according to \eqref{eq:giGMApprox} and  \eqref{eq:gijGMApprox} respectively) before the actual filter run. In that case, evaluation of these densities comes down to evaluating a precomputed Gaussian Mixture density.

\section{Conclusion}
In this paper we considered a Hidden Markov Model where the integrated continuous-time Markov chain can be observed at discrete time points perturbed by a Brownian motion. Using recent results on the asymmetric telegraph process, we derived exact formulas for a filter for the underlying continuous-time Markov chain, given this chain has only two different states. In case of more states we propose three approximate algorithms to compute the filter numerically. We investigated the performance of these filters with the help of two examples: one with two states and the other with five states. Problems like this arise for example when we have to filter the underlying economic state from observed asset prices.

\section*{Acknpwledgment}
We would like to thank our student Fehmi Mabrouk for some initial work on this topic during his Diploma Thesis which encouraged us to further pursue this topic.
 
\section*{Appendix}

{\em Proof of Theorem 3.1:}
The density of the integrated telegraph process $I_h$ is given in \cite{MKR94}:
$$ \Pop(I_h \in dx)/dx = e^{-\lambda h} \delta_h(x) + q_{h,\lambda}(x)$$ where $\delta_h(x)$ is the Dirac-measure in point $h$ and
$$ q_{h,\lambda}(x) = \frac\lambda 2 e^{-\lambda h} \Big[B_0\Big( \lambda \sqrt{(h^2-x^2)}\Big) + \Big(\frac{h+x}{h-x} \Big)^{\frac 12} B_1\Big( \lambda \sqrt{(h^2-x^2)}\Big)  \Big], \quad x\in [-h,h]$$
where the modified Bessel functions have been defined in (3.6). The density of $J_h$ given $\varepsilon_0=i$ is then by the density transformation formula
$$  \Pop_i(J_h \in dx)/dx = e^{-\lambda h} \delta_{h(a+b_i)}(x) + \frac{1}{|b_i|} q_{h,\lambda} \Big(\frac{x-ah}{b_i}\Big).$$
Finally the density of $Z_h=J_h+\sigma W_h$ given $\varepsilon_0=i$ is obtained by the convolution formula
\begin{eqnarray*}
  g_i(z,h) 
   &=& e^{-\lambda h} \phi_{\sigma\sqrt{h}}(z-h(b_i+a)) + \frac{1}{|b_i|} \int_{h (a-|b_i|)}^{h (a+|b_i|)}\phi_{\sigma\sqrt{h}}(z-x) q_{h,\lambda}\Big(\frac{x-ah}{b_i}\Big) dx.
\end{eqnarray*}
Thus, the first part of Theorem 3.1 is shown.
Next we determine the density of $I_h$ given $N_h$ is even. From this we can compute $g_{ii}$ as in the previous calculation. First note that
\begin{eqnarray*}
  \Pop( I_h \in dx | N_h \in 2 \N_0)/dx  &=& \frac{\Pop( I_h \in dx, N_h \in 2 \N_0)/dx}{\Pop(N_h \in 2 \N_0)} \\
   &=&  \frac{\sum_{k=0}^\infty \Pop( I_h \in dx, N_h = 2k)/dx}{\Pop(N_h \in 2 \N_0)} \\
   &=& \sum_{k=0}^\infty  \Pop( I_h \in dx | N_h =2k)/dx \frac{\Pop(N_h = 2k)}{\Pop(N_h \in 2 \N_0)}.
\end{eqnarray*}
For $k \ge 1$ the density of $I_h$ given $N_h=2k$ is
$$ \Pop(I_h \in dx | N_h =2k)/dx = \frac{(2k)!}{(k-1)!k!} (h^2-x^2)^{k-1} \frac{h+x}{(2t)^{2k}}, \quad x\in [-h,h],$$
see Lemma A.1 in \cite{MKR94}.
Moreover we have
\begin{eqnarray*}
  \Pop(N_h=2k)  &=& e^{-\lambda h} \frac{(\lambda h)^{2k}}{(2k)!} \\
    \Pop(N_h\in 2\N_0) &=& \sum_{k=0}^\infty   e^{-\lambda h} \frac{(\lambda h)^{2k}}{(2k)!} = e^{-\lambda h} \cosh(\lambda h).
\end{eqnarray*}
For $k=0$ the density of $I_h$ given $N_h=0$ has a point mass of one on $h$, i.e. $ \Pop(I_h \in dx | N_h =0)/dx=\delta_h(x)$ and $\Pop(N_h=0) = e^{-\lambda h}$.
Altogether we obtain:
\begin{eqnarray*}
  \Pop(I_h \in dx | N_h \in 2\N_0)/dx  &=& \delta_h(x) \cosh(\lambda h)^{-1} + \sum_{k=1}^\infty  \Pop(I_h \in dx | N_h =2k)/dx \frac{(\lambda h)^{2k}}{(2k)!\; \cosh(\lambda h)} \\
   &=& \cosh(\lambda h)^{-1} \delta_h(x) + q_{h,+}(x),
\end{eqnarray*}
where
$$ q_{h,+}(x)= \cosh(\lambda h)^{-1} \frac{\lambda}2 \Big(\frac{h+x}{h-x} \Big)^{\frac 12} B_1\Big(\sqrt{(h^2-x^2)} \lambda \Big), \quad x\in [-h,h]. $$

Next we determine the density $q_{h,-}(x)$ of $I_h$ given $N_h$ is odd.
For $k \ge 1$ the density of $I_h$ given $N_h=2k+1$ is
$$  \Pop(I_h \in dx | N_h =2k+1)/dx = \frac{(2k+1)!}{(k!)^2} \frac{(h^2-x^2)^k}{(2h)^{2k+1}}, \quad x\in [-h,h],$$
see Lemma A.1 in \cite{MKR94}. Moreover we have
\begin{eqnarray*}
  \Pop(N_h=2k+1)  &=& e^{-\lambda h} \frac{(\lambda h)^{2k+1}}{(2k+1)!} \\
    \Pop(N_h\in 2\N_0+1) &=& \sum_{k=0}^\infty   e^{-\lambda h} \frac{(\lambda h)^{2k+1}}{(2k+1)!} = e^{-\lambda h} \sinh(\lambda h).
\end{eqnarray*}
Altogether we obtain:
\begin{eqnarray*}
 q_{h,-}(x)  &=& \sum_{k=0}^\infty  \Pop(I_h \in dx | N_h =2k+1)/dx \frac{(\lambda h)^{2k+1}}{(2k+1)! \sinh(\lambda h)} \\
   &=& \sinh(\lambda h)^{-1}\frac{\lambda}2 B_0\Big(\sqrt{(h^2-x^2)} \lambda \Big), \quad x\in [-h,h].
\end{eqnarray*}
Hence the density of $J_h$ given $\varepsilon_0=\varepsilon_h=i$ is again obtained by the density transformation formula
$$ \Pop(J_h \in dx | \varepsilon_0=\varepsilon_h=i)/dx = \cosh(\lambda h)^{-1} \delta_{h(a+b_i)}(x) + \frac{1}{|b_i|} q_{h,+} \Big(\frac{x-ah}{b_i}\Big),$$
and finally the density of $Z_h=J_h +\sigma W_h$ given $\varepsilon_0=\varepsilon_h=i$ is
\begin{eqnarray*}
  g_{ii}(z,h) 
   &=& \phi_{\sigma\sqrt{h}}(z-h(b_i+a)) \cosh(\lambda h)^{-1} + \frac{1}{|b_i|} \int_{h (a-|b_i|)}^{h (a+|b_i|)}\phi_{\sigma\sqrt{h}}(z-x) q_{h,+}\Big(\frac{x-ah}{b_i}\Big) dx.
\end{eqnarray*}
Analogously we obtain the density of $J_h$ given $\varepsilon_0=i, \varepsilon_h=j$ with $i\neq j$ by
$$ \Pop(J_h \in dx | \varepsilon_0=i, \varepsilon_h=j)/dx = \frac{1}{|b_i|} q_{h,-} \Big(\frac{x-ah}{b_i}\Big)$$
and the density of $Z_h$ given $\varepsilon_0=i, \varepsilon_h=j$ is
\begin{eqnarray*}
  g_{ij}(z,h) &=&  \frac{1}{|b_i|} \int_{h (a-|b_i|)}^{h (a+|b_i|)}\phi_{\sigma\sqrt{h}}(z-x) q_{h,-}\Big(\frac{x-ah}{b_i}\Big) dx, \quad z\in\R.
\end{eqnarray*}

\vspace*{0.5cm}

{\em Proof of Theorem 3.2:}
The density of the asymmetric telegraph process $I_h$ given $\varepsilon_0=1$ is according to \cite{LR14}
$$ \Pop_1(I_h \in dx)/dx = e^{-\lambda_1 h} \delta_h(x) + p_{h}(x)$$ where $\delta_h(x)$ is the Dirac-measure in point $h$ and
\begin{eqnarray*} p_{h}(x) &=& \frac{\lambda_1}2 \exp\Big(-\lambda_1\big(\frac{x+h}2) -\lambda_2\big(\frac{h-x}2\big) \Big) \\
&& \cdot \Big[B_0\Big( \sqrt{\lambda_1\lambda_2 (h^2-x^2)}\Big) + \Big(\frac{(h+x)\lambda_2}{(h-x)\lambda_1} \Big)^{\frac 12} B_1\Big( \sqrt{ \lambda_1\lambda_2 (h^2-x^2)}\Big)  \Big], \quad x\in [-h,h].
\end{eqnarray*}
From this, the density of $Z_h$ given $\varepsilon_0=1$ follows as in the previous proof.

Next we have to determine the density of $Z_h$ given $\varepsilon_0=\varepsilon_h=1$. Note that we first have:
\begin{eqnarray*}
  g_{11}(z,h) &=& \Pop( Z_h \in dx | \varepsilon_0=\varepsilon_h=1)/dx = \frac{\Pop( Z_h \in dx,  \varepsilon_h=1 | \varepsilon_0=1)/dx}{\Pop(\varepsilon_h=1 |\varepsilon_0=1)}.
\end{eqnarray*}
From (3.7) we obtain that
$$ \Pop(\varepsilon_h=1 |\varepsilon_0=1) = p_{11}(h) = \frac{\lambda_2}{\lambda_1+\lambda_2}+\frac{\lambda_1}{\lambda_1+\lambda_2} e^{-(\lambda_1+\lambda_2)t}.$$
From the results in \cite{LR14} we can derive $ \Pop( I_h \in dx,  \varepsilon_h=1 | \varepsilon_0=1)/dx$. We have that
$$ \Pop( I_h \in dx,  \varepsilon_h=1 | \varepsilon_0=1)/dx = e^{-\lambda_1 h} \delta_h(x) + q_{h,+,1}(x)$$ where
\begin{eqnarray*} q_{h,+,1}(x) &=& \frac{1}2 \exp\Big(-\lambda_1\big(\frac{x+h}2) -\lambda_2\big(\frac{h-x}2\big) \Big) \\
&& \cdot \Big[ \sqrt{\lambda_1\lambda_2}\Big(\frac{h+x}{h-x} \Big)^{\frac 12} B_1\Big( \sqrt{ \lambda_1\lambda_2 (h^2-x^2)}\Big)  \Big], \quad x\in [-h,h].
\end{eqnarray*}
The density $g_{11}(z,h)$ follows by joining the results.\\

We proceed similar for the density of $Z_h$ given $\varepsilon_0=1, \varepsilon_h=2$. Here we have
\begin{eqnarray*}
  g_{12}(z,h) &=& \Pop( Z_h \in dx | \varepsilon_0=1, \varepsilon_h=2)/dx = \frac{\Pop( Z_h \in dx,  \varepsilon_h=2 | \varepsilon_0=1)/dx}{\Pop(\varepsilon_h=2 |\varepsilon_0=1)}.
\end{eqnarray*}
From (3.7) we obtain
$$ \Pop(\varepsilon_h=2 |\varepsilon_0=1) = p_{12}(h) = \frac{\lambda_1}{\lambda_1+\lambda_2}-\frac{\lambda_1}{\lambda_1+\lambda_2} e^{-(\lambda_1+\lambda_2)t}$$
and for the density of $I_h$ and $\varepsilon_h=2$ given $\varepsilon_0=1$ in the asymmetric case note that
\begin{eqnarray*} \Pop( I_h \in dx,  \varepsilon_h=2 | \varepsilon_0=1)/dx  &=& \frac{\lambda_1}2 \exp\Big(-\lambda_1\big(\frac{x+h}2) -\lambda_2\big(\frac{h-x}2\big) \Big) \\
&& \cdot B_0\Big( \sqrt{\lambda_1\lambda_2 (h^2-x^2)}\Big) , \quad x\in [-h,h].
\end{eqnarray*}
Using the usual transformation yields the statement.

\bibliography{Bibfilter,ig-literature}
\bibliographystyle{mybibstylefile1}
\end{document}